\documentclass[onefignum,onetabnum]{siamart190516}

\usepackage[utf8]{inputenc}
\usepackage{amsmath}

\usepackage{amsfonts}
\usepackage{amssymb,mathtools}
\usepackage{latexsym}
\usepackage{bbm}
\usepackage{listings}
\usepackage{graphicx,xcolor}
\usepackage{subfig,caption,booktabs}
\usepackage{multirow}
\usepackage{enumitem}
\usepackage[toc,page]{appendix}
\usepackage[algo2e,algoruled]{algorithm2e}
\usepackage[sort]{cite}
\usepackage{siunitx}
\newcommand{\rd}{\,\mathrm{d}}
\DeclareMathOperator{\Tr}{Tr}
\DeclareMathOperator{\Hessian}{Hess}
\DeclareMathOperator{\ReLU}{ReLU}
\newcommand{\forcond}{$\ell=1$ \KwTo the number of iterations}

\newsiamthm{exmp}{Example}
\newsiamremark{remark}{Remark}

\DeclareMathOperator{\dist}{dist}

\headers{Actor-Critic Method for High Dimensional Static HJB equations}{Mo Zhou, Jiequn Han, and Jianfeng Lu}

\begin{document}

\title{Actor-Critic Method for High Dimensional Static Hamilton--Jacobi--Bellman Partial Differential Equations Based on Neural Networks \thanks{Date: February 20, 2021. Please address correspondence to Jiequn Han or Jianfeng Lu. \funding{The work of JL and MZ is supported in part by National Science Foundation via grant DMS-2012286. The authors are grateful for computing time at the Terascale Infrastructure for Groundbreaking Research in Science and Engineering (TIGRESS) of Princeton University.
}}}

\author{Mo Zhou\thanks{Department of Mathematics, Duke University (\email{mo.zhou366@duke.edu}).} 
\and Jiequn Han\thanks{Department of Mathematics, Princeton University (\email{jiequnhan@gmail.com}).}
\and Jianfeng Lu\thanks{Department of Mathematics, Department of Physics, and Department of Chemistry, Duke University (\email{jianfeng@math.duke.edu}).}}

\maketitle 

\begin{abstract}
We propose a novel numerical method for high dimensional Hamilton--Jacobi--Bellman (HJB) type elliptic partial differential equations (PDEs). The HJB PDEs, reformulated as optimal control problems, are tackled by the actor-critic framework inspired by reinforcement learning, based on neural network parametrization of the value and control functions. Within the actor-critic framework, we employ a policy gradient approach to improve the control, while for the value function, we derive a variance reduced least-squares temporal difference method using stochastic calculus. To numerically discretize the stochastic control problem, we employ an adaptive step size scheme to improve the accuracy near the domain boundary. Numerical examples up to $20$ spatial dimensions including the linear quadratic regulators, the stochastic Van der Pol oscillators, the diffusive Eikonal equations, and fully nonlinear elliptic PDEs derived from a regulator problem are presented to validate the effectiveness of our proposed method.
\end{abstract}


\section{Introduction}

The Hamilton-Jacobi-Bellman (HJB) equation is an important family of partial differential equations (PDEs), given its connection with optimal control problems that lead to a wide range of applications.
The unknown in the HJB equation can be viewed as the total expected value function for optimal control problems. The equation can be derived from the dynamic programming principle pioneered by Bellman~\cite{bellman1966dynamic}, which gives a necessary and sufficient condition of the optimality.
Theoretical results for the existence and uniqueness of the HJB equations are well established; see, e.g., \cite{yong1999stochastic}. 
From the viewpoint of stochastic control, the relationship between the viscosity solution of the HJB equations and the backward stochastic differential equations (BSDEs) is introduced in \cite{crandall1992user,pardoux1990adapted,peng1991probabilistic,pardoux1992backward,pardoux1998backward}. 

The wide applications of HJB equations call for efficient numerical algorithms. Various numerical approaches have been developed in the literature, including the monotone approximation scheme \cite{barles2002convergence, forsyth2007numerical}, the finite volume method \cite{wang2003numerical,richardson2006numerical}, and the Galerkin method \cite{beard1997galerkin, beard1998approximate}.
In \cite{osher1988fronts}, nonoscillatory schemes are developed to solve the HJB equations exploring the connection with hyperbolic conservation laws.
The HJB equations related to reachability problems are studied in \cite{mitchell2003overapproximating, mitchell2005time, lygeros2004reachability}.
A general survey for classical methods to solve the optimal control problem numerically can be found, e.g., in \cite{rao2009survey}.
While these conventional approaches have been quite successful, they fall short for solving HJB equations in high dimensions due to the curse of dimensionality~\cite{bellman1966dynamic}: the computational cost goes up exponentially with the dimensionality.
Many works attempt to mitigate this fundamental difficulty by leveraging dimension reduction techniques such as proper orthogonal decomposition, sparse grid, pseudospectral collocation, and tensor decomposition~(see e.g.,~\cite{kunisch2004hjb,kang2017mitigating,kalise2018polynomial,dolgov2019tensor,oster2019approximating}).
The performance of these algorithms heavily depends on how well the low dimensional representation matches the solutions, and is typically problem dependent and thus with limited applicability.

To better address the challenge of high dimensionality, a promising direction is to consider the artificial neural network as a more flexible and efficient function approximation tool. This topic has received a considerable amount of attention and been a rapidly developing field in recent years. Several numerical approaches for high dimensional PDEs based on neural network parametrization have been proposed; see e.g., the reviews \cite{han2020algorithms, beck2020overview} and references therein. 

For HJB type equations and related optimal control problems, the most tightly connected approach to our work is the deep BSDE method \cite{han2018solving, weinan2017deep}, which reformulates parabolic PDEs as control problems using BSDEs, and uses deep neural network parametrization for the solution and control to solve this problem. Theoretical results for convergence of this method are studied in \cite{han2020convergenceBSDE}. 
The deep BSDE method and its variants have been applied to solve HJB type equations, stochastic control problems, and differential games (see e.g., \cite{han2018solving, weinan2017deep,chan2019machine, henry2017deep, pereira2019learning, pham2021neural,nusken2020solving,kremsner2020deep,ji2020three,han2020deep,han2020convergence}).
Numerical algorithms for solving high dimensional deterministic and stochastic control problems based on other forms combined with deep learning approximation have also been investigated in~\cite{han2016deep,nakamura2019adaptive,becker2019deep,han2021recurrent}.

While some methods mentioned above have been successful in solving PDEs in high dimensions, there have been two issues that remain to be addressed. On the one hand, most of these works concern parabolic PDEs of finite time horizon (often of order one), while only a few works investigate the static elliptic HJB equations corresponding to control problems with infinite time horizon. On the other hand, most existing works consider equations where the optimal controls are explicitly known given the value function or without controls, while there are many important HJB type equations for which the optimal control is cast through an optimization problem and hence implicit. 
Recently, an algorithm for a high dimensional finite-time horizon stochastic control problem with implicit optimal control is considered in~\cite{ji2020deep}, based on the deep BSDE formulation associated with the stochastic maximum principle.
In this paper, we take a different approach and focus on solving the static elliptic type HJB equation with implicit control, in which the above two challenges are compounded.

Our proposed numerical method is heavily inspired by the literature on reinforcement learning (RL)~\cite{sutton2018reinforcement}, which is of course closely related to control problems. Our motivation for borrowing techniques from RL is due to the impressive revolution and great success in recent years in deep RL by utilizing neural network parametrization~\cite{mnih2013playing,silver2016mastering,duan2016benchmarking}. 
In the RL context, the control problem is usually formulated as a Markov decision process (MDP) on discrete time and state space. If the model is given, finding the optimal policy can be viewed as solving a discrete HJB equation. It is then natural to ask whether algorithms developed in the RL context can be generalized to the context of solving high dimensional HJB equations.

In this paper, we reformulate the HJB type fully nonlinear elliptic PDEs into stochastic control problems and leverage the actor-critic framework in conjunction with a neural network approximation to solve the equations.
The actor-critic methods are a class of algorithms in RL \cite{sutton2018reinforcement}. These algorithms iteratively evaluate and improve the current policies (i.e., controls) until final convergence.
The critic refers to the value function of a given policy. The process of estimating the critic is called policy evaluation.
The most common algorithms for policy evaluation are temporal difference (TD) methods \cite{konda2000actor,bhatnagar2009natural, vamvoudakis2010online}, or their variants, such as the TD$(\lambda)$ \cite{degris2012off} and the least-squares TD (LSTD) \cite{peters2008natural, maei2010toward, boyan1999least}.
The actor refers to the policy function, and we need to make policy improvement based on a given value function.
In this case, the most popular method is policy gradient \cite{konda2000actor,bhatnagar2006reinforcement,abdulla2007parametrized, degris2012off,wang2016sample} and their variants, such as natural policy gradient \cite{peters2008natural,bhatnagar2009natural}.
In this work, we propose a variance reduced version of the LSTD method for policy evaluation derived using stochastic calculus. We also adapt the policy gradient method for policy improvement to the continuous-time stochastic control problem.

The rest of this paper is organized as follows. In Section 2, we provide a theoretical background for the optimal control problems and formulate the actor-critic framework for continuous-time stochastic control. In Section 3, we introduce the numerical algorithm to solve the optimal control problem. Numerical examples are presented in Section 4. We conclude in Section 5 with an outlook for future works.

\section{Theoretical background of the actor-critic framework}
\subsection{Control formulation of elliptic equations} \label{sec:control_theory}

Consider the following fully nonlinear elliptic PDE
\begin{equation}\label{eqn:NonlinearPDE}
\inf_{u \in U} \left[\frac{1}{2} \Tr \left( \sigma \sigma^{\top} \Hessian(V) \right) (x, u) +  b(x,u)^{\top} \nabla V(x) + f(x,u)\right]  - \gamma V(x) = 0  \quad \text{in }  \Omega,
\end{equation}
with boundary condition $V(x) = g(x)$ on $\partial \Omega$. Here the state space $\Omega$ is an open, connected set in $\mathbb{R}^d$ with piecewise smooth boundary, and the control space $U$ is a convex closed domain in $\mathbb{R}^{d_u}$. We assume that $V(x) \in C^2(\overline{\Omega})$, $f(x,u) \in C(\overline{\Omega} \times U)$, $b(x,u) \in C(\overline{\Omega} \times U; \mathbb{R}^d)$,
$\sigma(x, u) \in C(\overline{\Omega} \times U; \mathbb{R}^{d \times d_w})$
with 
$\sigma(x, u) \sigma^{\top} (x, u)$
being uniformly elliptic and bounded, and $\gamma \geq 0$ is a constant. Here and in the following, we use $\nabla$ and $\Hessian$ to denote the gradient and Hessian operators.

As a starting point of our approach, we reformulate the above elliptic equation as an optimal control problem. Let $(\widetilde\Omega, \mathcal{F}, \{\mathcal{F}_t\}_{t \ge 0}, \textbf{P})$ be a filtered probability space. Consider the following stochastic differential equation (SDE)
\begin{equation} \label{eqn:SDE}
\rd X_t = b(X_t, u_t) \rd t + \sigma (X_t, u_t) \rd W_t 
\end{equation}
with initial condition $X_0 = x\in \Omega$, where $u_t \in U \subset \mathbb{R}^{d_u}$ is an $\mathcal{F}_t$-adapted control field and $W_t$ is a $d_w$-dimensional $\mathcal{F}_t$-standard Brownian motion.
As we solve the equation in the domain $\Omega$, we define a stopping time 
\begin{equation} \label{eqn:stopping time}
\tau = \inf\{ t: X_t \notin \Omega \}.
\end{equation}
It is a standard result that $\tau < \infty$ a.s.; see, for example, \cite{klebaner2005introduction}.

We then consider an optimal control problem to minimize the following cost functional
\begin{equation} \label{eqn:cost_functional}
J^u(x) = \mathbb{E} \Bigl[\int_0^{\tau} f(X_s, u_s) e^{-\gamma s} \rd s + e^{-\gamma \tau} g(X_{\tau}) \mid X_0^u = x  \Bigr].
\end{equation}
In this cost functional, $f$ can be interpreted as running cost, $g$ is the terminal cost when the SDE hits the boundary $\partial \Omega$, and $\gamma$ is the discount rate. 

The control $u$ is chosen over the set of stochastic processes that have values in $U$ and are adapted to the filtration $\mathcal{F}_t$.
Define
\begin{equation}
V(x) = \inf_{u} J^u(x)
\end{equation}
as the optimal value function (i.e., optimal cost-to-go function). According to standard results in stochastic control theory~\cite{yong1999stochastic}, $V$ satisfies the time-independent HJB equation 
\begin{equation} \label{eqn:HJB}
\inf_{u}  \bigl\{ \mathcal{L}^u V(x, u) + f(x,u) - \gamma V(x) \bigr\} =0
\end{equation}
in $\Omega$ with boundary condition $V(x) = g(x)$ on $\partial \Omega$, where 
\begin{equation*} 
\mathcal{L}^u V(x) = \frac{1}{2} \Tr \left( \sigma \sigma^{\top} \Hessian(V) \right) (x, u) + b(x,u)^{\top} \nabla V(x)
\end{equation*}
is the generator of the controlled SDE \eqref{eqn:SDE}. 
Note that the HJB equation \eqref{eqn:HJB} coincides with the original PDE \eqref{eqn:NonlinearPDE} and, hence, we can solve the PDE \eqref{eqn:NonlinearPDE} by solving the optimal control problem to obtain the optimal value function. 

\subsection{Actor-critic method in stochastic optimal control problem} \label{sec:diecrete_RL}

Our approach for solving the optimal control problem is based on the actor-critic framework. In such methods, one solves for both the value function and control field. The control (i.e., policy in the RL terminology) is known as the \emph{actor}, while the value function corresponding to the control is known as the \emph{critic} since it is used to evaluate the optimality of the control. Accordingly, the actor-critic algorithms consist of two parts: policy evaluation for the critic and policy improvement for the actor. 
While many approaches have been developed under the actor-critic framework~\cite{konda2000actor,abdulla2007parametrized,peters2008natural,bhatnagar2009natural,vamvoudakis2010online,degris2012off,wang2016sample}, we will focus on simple and perhaps the most popular algorithms: TD learning for the value function given a policy and policy gradient for improving the control. 

\subsubsection{TD for discrete Markov decision processes} 

To better convey the idea, let us first briefly recall the algorithms for the discrete-time MDP with finite state and action space; more details can be found in e.g.,~\cite{sutton2018reinforcement}.
The MDP starts with some initial state $S_0$ in the state space $\mathcal{S}$, possibly sampled according to a distribution. 
At time $t \in \mathbb{N}$, given the current state $S_t$, the agent picks an action $A_t$ in the action set $\mathcal{A}$ according to a policy. We assume that the policy is deterministic, i.e., the policy is a map $\pi$ from the state space $\mathcal{S}$ to the action space $\mathcal{A}$:
\begin{equation} \label{eqn:MDPaction}
    A_t = \pi(S_{t}). 
\end{equation}
After the action $A_t$ is chosen, the system state will transit to $S_{t+1}$, according to a probability transition function 
\begin{equation}
    \mathbb{P}( S_{t+1} = s' \mid S_t = s, A_t = a) = p(s' \mid s, a). 
\end{equation}
The action also incurs a cost $R_{t+1}$, which we assume to be given by a deterministic function of the previous state $S_t$, action $A_t$, and the current state $S_{t+1}$:
\begin{equation} \label{eqn:MDPsingle_cost}
R_{t+1} = f(S_t, A_t, S_{t+1}).
\end{equation}
The goal of the MDP problem is to choose the best policy to minimize the expected total discounted cost 
\begin{equation} \label{eqn:MDPcost}
\mathbb{E}_{S_0 \sim \mu, \pi} \, \Bigl[\sum_{t=1}^{\infty} \beta^{t-1} R_t\Bigr],
\end{equation}
where $\beta \in (0,1)$ is a discount factor, $\mu$ is the distribution of the initial state $S_0$, and we have used $\mathbb{E}_{\pi}$ to indicate the dependence on the transition dynamics on the choice of the policy $\pi$. 

\smallskip 

To solve the MDP problem, it is convenient to introduce the (state) value function 
w.r.t. a policy $\pi$ as the expected cost starting at states $s$ under that policy:
\begin{equation}\label{eqn:MDPvalue}
V^{\pi}(s) = \mathbb{E}_{\pi}\,\Bigl[\sum_{t=1}^{\infty} \beta^{t-1} R_t \mid S_0 = s\Bigr].
\end{equation}
By the dynamic programming principle~\cite{bellman1966dynamic}, for any given policy $\pi$, the value function satisfies
\begin{equation}\label{eqn:MDPvalue2}
V^{\pi}(s) = \mathbb{E}_{\pi}\,\Bigl[\sum_{t=1}^{n} \beta^{t-1} R_{t} + \beta^{n} V^{\pi}(S_{n}) \mid S_0 = s\Bigr]
\end{equation}
for any $n \geq 1$.
In order to minimize the total cost \eqref{eqn:MDPcost}, we search for an optimal policy $\pi^*$ that satisfies for all $\pi$,
\begin{equation}
    V^{\pi^*}(s) \le V^{\pi}(s) \qquad  \forall\, s \in \mathcal{S}. 
\end{equation}
Specifically, by the optimality principle, we have  
\begin{equation}\label{eqn:MDP_DP}
V^{\pi^*}(s) = \min_{\{a_t\} \subset \mathcal{A}} \mathbb{E}\,\Bigl[\sum_{t=1}^{n} \beta^{t-1} R_{t} + \beta^{n} V^{\pi^*}(S_{n}) \mid S_0 = s, A_{t} = a_t, t=0,1,\cdots,n-1].
\end{equation}
Note that while in \eqref{eqn:MDPvalue2} and \eqref{eqn:MDP_DP} the right-hand side starts at time $0$, we can start at any time and run the process for $n$ steps due to stationarity. 

Let us make a couple of remarks for the setup of the discrete MDP used here. 
First, in RL, reward is usually used instead of cost and, hence, one maximizes the total reward instead of minimizing the cost; evidently, the two viewpoints are equivalent up to a change of sign. We use cost, which is more in line with the control literature and also our problem in the continuous setting. 
Second, the cost $R_t$ is not necessarily a deterministic function as in \eqref{eqn:MDPsingle_cost}, but may follow some probability distribution together with the next state:
\begin{equation}
     \mathbb{P}( S_{t+1} = s', R_{t+1}=r \mid S_t = s, A_t = a) = p(s', r \mid s, a). 
\end{equation}
Moreover, the policy can also be probabilistic rather than deterministic as assumed in \eqref{eqn:MDPaction}. 
We choose the simplified setting for the cost and policy to make it consistent with our continuous optimal control setting.
Finally, we use an MDP without stopping time and thus without terminal cost for simplicity. The adaptation to our PDE setup will be discussed below in Sections~\ref{sec:contTD} and \ref{sec:contPG}.

\medskip 

TD learning is a class of algorithms that evaluate a given control (i.e., policy) by updating the value function, combining a Monte Carlo estimate of the running cost over a time period and the dynamic programming principle for the future cost-to-go. 
For a policy $\pi$ to be evaluated, with a given trajectory $\{S_t, t\geq0\}$, we update the value function at each $t$ by 
\begin{equation}\label{eqn:MDP_TD_update}
\widehat{V}^{\pi}(S_t) \longleftarrow \widehat{V}^{\pi}(S_t) + \alpha \left( \sum_{k=1}^{n} \beta^{k-1} R_{t+k} + \beta^{n} \widehat{V}^{\pi}(S_{t+n}) - \widehat{V}^{\pi}(S_t) \right),
\end{equation}
where $\alpha$ is the learning rate and $\widehat{V}^{\pi}$ on the right-hand side is the current estimate of the value function. 
In \eqref{eqn:MDP_TD_update}, we only update the value function at the state $S_t$ and the value at other states remain unchanged. 
In practice, this update of the value function is usually done for multiple trajectories.

In the above updating rule, we have used the $n$-step TD  $\mathrm{TD}_n^{\pi}(S_t)$, defined as
\begin{equation}\label{eqn:MDP_TD}
\mathrm{TD}_n^{\pi}(S_t) = \sum_{k=1}^{n} \beta^{k-1} R_{t+k} + \beta^{n} \widehat{V}^{\pi}(S_{t+n}) - \widehat{V}^{\pi}(S_t),
\end{equation}
which depends on the trajectory of length $n+1$ ($t$ to $t+n$) from the starting state $S_t$.
$\mathrm{TD}_n^{\pi}$ can be understood as an indicator of the inconsistency between the current estimate of the value function with a sampled value using $n$-steps of the MDP, since according to \eqref{eqn:MDPvalue2}, $\mathbb{E}_{\pi} \mathrm{TD}_n^{\pi}$ vanishes if $\widehat{V}^\pi$ agrees with the true value function $V^{\pi}$. Hence, TD learning can be viewed as a stochastic fixed point iteration for the value function. 

When function approximation is used for the value function, in particular nonlinear approximations such as neural networks, an alternative approach, the LSTD is often used to overcome potential divergence problems of TD learning \cite{peters2008natural, boyan1999least}. Instead of the stochastic fixed point updating formula as  \eqref{eqn:MDP_TD_update}, in the LSTD method, the parameters are optimized to minimize the squares of the TD error as a loss function. More specifically, if the value function is parametrized as $V^{\pi}(\cdot; \theta_V)$, we solve for $\theta_V$ by 
\begin{equation}\label{eq:optLSTD}
    \min_{\theta_V} \mathbb{E}_{S_0 \sim \mu, \pi} \Bigl[ \Bigl( \sum_{t=1}^{n} \beta^{t-1} R_{t} + \beta^{n} V^{\pi}(S_{n}; \theta_V) - V^{\pi}(S_0; \theta_V) \Bigr)^2 \mid S_0  \Bigr],
\end{equation}
where $\mu$ is some initial distribution for the state $S_0$.
In practice, \eqref{eq:optLSTD} is often solved using the stochastic gradient descent method. 
Such method has been proved successful in e.g., \cite{mnih2013playing,duan2016benchmarking}. 

\subsubsection{Policy gradient for discrete MDP}\label{sec:MDP_PG}

Policy gradient is a class of methods to learn parametrized policies through gradient based algorithms. Assume we consider a class of (deterministic) policy parametrized as 
\begin{equation}
A_t(S_t) = \pi(S_t; \theta_{\pi}),
\end{equation}
where $\theta_\pi$ denotes a collection of parameters and $\pi(\cdot; \theta)$ is a chosen nonlinear parametrization. 

To find an optimal policy, we aim to minimize the objective function (cf.~\eqref{eqn:MDP_DP})
\begin{equation}\label{eq:MDP_actor}
J(\theta_{\pi}) = \mathbb{E}_{S_0\sim \mu, \pi(\cdot; \theta_{\pi})} \bigl[ \sum_{k=1}^{n} \beta^{k-1} R_{k} + \beta^{n} \widehat{V}^{\pi}(S_{n}) \bigr]
\end{equation}
w.r.t. the collective parameter $\theta_\pi$. Note that \eqref{eq:MDP_actor} explicitly takes into account the cost of the first $n$ steps, while using an (approximate) value function $\widehat{V}^{\pi}$ for the future cost after $n$ steps, coming from e.g., the TD learning algorithm. The objective function \eqref{eq:MDP_actor} can be optimized using stochastic gradient method. Using a stochastic estimate of the gradient $ \widehat{\nabla J} \approx \nabla_{\theta_{\pi}} J$, so that the parameter is updated as  
\begin{equation}\label{eqn:MDP_policy_update}
\theta_{\pi} \gets \theta_{\pi} - \alpha \widehat{\nabla J}
\end{equation}
with suitable learning rate $\alpha$. Note that in principle, one needs to differentiate all terms involved in \eqref{eq:MDP_actor} w.r.t. $\theta_{\pi}$; however, in practice, in the actor-critic framework, one typically leaves out the derivative of $\widehat{V}^{\pi}$ w.r.t. $\pi$, as it is impractical to compute since $\widehat{V}^{\pi}$ is obtained using e.g., TD learning. Nevertheless, we would still need to differentiate $\widehat{V}^{\pi}(S_{n})$ w.r.t. $S_n$, as the state $S_{n}$ is affected by the choice of the policy, thus
\begin{equation*}
    \frac{\partial \widehat{V}^{\pi}(S_n)}{\partial \theta_{\pi}} \stackrel{\cdot}{=} \frac{\partial \widehat{V}^{\pi}}{\partial S_n} \frac{\partial S_n}{\partial \theta_{\pi}},
\end{equation*}
where $\stackrel{\cdot}{=}$ indicates that the (functional) derivative $\frac{\delta \widehat{V}^{\pi}}{\delta \pi}$ is omitted. 
Dropping this term is often applied in actor-critic algorithms. Some justifications can be found in~\cite{degris2012off}: Under certain conditions, the approximated gradient is still in the direction of improving the performance and the set of critical points of the objective function coincides with the set of zero approximated gradients.

Since we consider deterministic policies, the policy gradient approach discussed above is in the same spirit as the \emph{deterministic policy gradient algorithm} proposed in \cite{silver2014deterministic}. One difference is that we use the state value function $V(s)$ while \cite{silver2014deterministic} uses the state-action value function ($Q$-function) $Q(s,a)$. Another difference is that the objective function used in \cite{silver2014deterministic}  is based on the stationary distribution of a state-action pair given the policy, while we roll out a trajectory for the cost function (combined with using an estimated value function for future cost), which is more suitable to an actor-critic framework. Our approach is also easier to generalize to the continuous setting, which will be discussed in Section~\ref{sec:contPG}.

\subsubsection{TD for continuous optimal control problems} \label{sec:contTD}

We now introduce how to adapt the above algorithmic ideas to the continuous setting.

Given a control function $u(x) \in C(\overline{\Omega})$ (which corresponds to $\pi$ in the discrete setting), the corresponding value function is given by
\begin{equation} \label{eqn:value}
V^u(x) = \mathbb{E}_u\, \Bigl[\int_0^{\tau} f(X_s, u(X_s)) e^{-\gamma s} \rd s + e^{-\gamma \tau} g(X_{\tau}) \mid X_0 = x  \Bigr]. 
\end{equation}
Here $\mathbb{E}_u$ indicates expectation w.r.t. the trajectory (with a fixed policy $u$).
This is just the cost functional in \eqref{eqn:cost_functional} with a specific control policy. In the continuous setting, the dynamical programming principle indicates that the value function $V^u$ satisfies the PDE (see e.g.,~\cite{yong1999stochastic})
\begin{equation}\label{eqn:linearPDE}
\frac{1}{2} \Tr \left( \sigma \sigma^{\top} \Hessian(V^u) \right)(x, u(x)) + b(x,u(x))^{\top} \nabla V^u(x) + f(x,u(x)) - \gamma V^u(x) = 0 \, \text{in }\Omega
\end{equation}
with boundary condition $V^u(x) = g(x)$ on $\partial \Omega$.

To better convey the idea, we first consider a fixed time interval $[0, T]$ with $T > 0$ and neglect the stopping time and also the domain boundary. Necessary modifications regarding the stopping time and boundary will be explained below. 
Applying It\^o's formula to $e^{-\gamma t} V^u(X_t)$, we get 
\begin{multline} \label{eqn:Ito1}
e^{-\gamma T} V^u(X_T)  =  V^u(X_0) + \int_0^T e^{-\gamma s} \Bigl[\frac{1}{2} \Tr \bigl( \sigma \sigma^{\top} \Hessian(V^u) \bigr)(X_s, u(X_s)) \\
\hfill + b(X_s,u(X_s))^{\top} \nabla V^u(X_s)
- \gamma V^u(X_s)\Bigr] \rd s \\
+ \int_0^T e^{-\gamma s} \nabla V^u(X_s)^{\top} \sigma(X_s, u(X_s)) \rd W_s.
\end{multline}
Combined with the PDE \eqref{eqn:linearPDE}, \eqref{eqn:Ito1} gives
\begin{multline} \label{eqn:Ito2}
V^u(X_0) = \int_0^T e^{-\gamma s} f(X_s, u(X_s)) \rd s\\ - \int_0^T e^{-\gamma s} \nabla V^u(X_s)^{\top} \sigma(X_s, u(X_s)) \rd W_s + e^{-\gamma T} V^u(X_T).
\end{multline}
The term $\int_0^T e^{-\gamma s} \nabla V^u(X_s)^{\top} \sigma(X_s, u(X_s)) \rd W_s$ is a martingale w.r.t. $T$ because $\nabla V$ and $\sigma$ are bounded according to our assumptions \cite{klebaner2005introduction}. Therefore, taking the expectation, we arrive at
\begin{equation} \label{eqn:cond_exp1}
V^u(X_0) = \mathbb{E}_u\; \Bigl[ \int_0^T e^{-\gamma s} f(X_s, u(X_s)) \rd s + e^{-\gamma T} V^u(X_T) \mid X_0 \Bigr].
\end{equation}
We observe that this is the analog of \eqref{eqn:MDPvalue2} in the continuous setting, where the unit time discount $e^{-\gamma}$ is the analog of the discount factor $\beta$ in the discrete setting. 
Compared with the discrete time setting, besides \eqref{eqn:cond_exp1}, we have in addition \eqref{eqn:Ito2} before taking expectation, thanks to It\^o's lemma.
Exploiting the two identities, analogously to the discrete case, we define two versions of TD in the continuous setting as
\begin{align}
& \mathrm{TD}_{1}^u  = \int_0^T e^{-\gamma s} f(X_s,u(X_s)) \rd s - \int_0^T e^{-\gamma s} \nabla V(X_s)^{\top} \sigma(X_s,u(X_s)) \rd W_s \\
& \hspace{10em} + e^{-\gamma  T}  V(X_T) - V(X_0), \nonumber \\
& \mathrm{TD}_2^u = \int_0^T e^{-\gamma s} f(X_s,u(X_s)) \rd s + e^{-\gamma  T}  V(X_T) - V(X_0).
\end{align}
Note that both $\mathrm{TD}_1$ and $\mathrm{TD}_2$ depend on the trajectory $X_t$; in particular, they should be viewed as random variables, while we suppress such dependence in the notation. 
From \eqref{eqn:Ito2} and \eqref{eqn:cond_exp1}, if $V$ is the exact value function corresponding to the control $u$, we have 
\begin{align}
    & \mathrm{TD}_1^u = 0,\quad \mathbb{P}\text{-a.s.} \\
    & \mathbb{E}_u \, \mathrm{TD}_2^u = 0.
\end{align}
Note in particular that $\mathrm{TD}^u_1$ vanishes without taking the expectation for the exact value function while $\text{Var}(\mathrm{TD}^u_2) = \mathbb{E}_u [\int_0^T e^{-2 \gamma s} |\nabla V(X_s)^{\top} \sigma(X_s,u(X_s))|^2 \rd s]  > 0$ if $\nabla V^{\top}(x) \sigma(x,u)\not\equiv 0$. Moreover, as the difference between $\mathrm{TD}_1$ and $\mathrm{TD}_2$ is given by a martingale term, for any approximate value function, we have 
\begin{equation*}
    \mathbb{E}_u \mathrm{TD}_1^u = \mathbb{E}_u \mathrm{TD}_2^u. 
\end{equation*}
Now let us introduce two loss functionals for the critic in the spirit of LSTD:
\begin{align}
L_1(V) &= \mathbb{E}_{X_0 \sim \mu, u} \, \bigl( \mathrm{TD}_1^u \bigr)^2  \label{eqn:Critic1} \\
& = \mathbb{E}_{X_0 \sim \mu, u} \, \biggl[\Bigl( \int_0^{T \land \tau} e^{-\gamma s} f(X_s,u(X_s)) \rd s  \notag\\
& \hspace{1em} - \int_0^{T \land \tau} e^{-\gamma s} \nabla V(X_s)^{\top} \sigma(X_s,u(X_s)) \rd W_s + e^{-\gamma (T \land \tau)}  V(X_{T \land \tau}) - V(X_0) \Bigr)^2\biggr],\notag   \\
L_2(V) &= \mathbb{E}_{X_0 \sim \mu, u}\, \bigl( \mathrm{TD}_2^u \bigr)^2 \label{eqn:Critic2}
\\
& = \mathbb{E}_{X_0 \sim \mu, u} \,  \biggl[\Bigl(  \int_0^{T \land \tau} e^{-\gamma s} f(X_s,u(X_s)) \rd s + e^{-\gamma (T \land \tau)}  V(X_{T \land \tau}) - V(X_0) \Bigr)^2\biggr], \notag
\end{align}
where $\mu$ is some initial distribution for $X_0$ and we have also taken into account the stopping time $\tau$ when the process hits the domain boundary. Here the two losses are viewed as functionals of the value function $V$, the finite-dimensional function approximation will be discussed in the next section.

\smallskip 

The stochastic gradient method is used to minimize the loss function in LSTD to find the best approximation of the value function. Written in terms of functional variations, this amounts to approximating
\begin{align}
    & \mathbb{E}_{X_0 \sim \mu, u} \frac{\delta \bigl(\mathrm{TD}_1^u)^2 }{\delta V} \approx \frac{\delta \bigl(\mathrm{TD}_1^u)^2 }{\delta V}\bigl(X_t\bigr) = 2  \mathrm{TD}_1^u\bigl(X_t\bigr) \frac{\delta \mathrm{TD}_1^u}{\delta V} \bigl( X_t \bigr), \label{eq:TDgrad1} \\  
    & \mathbb{E}_{X_0 \sim \mu, u} \frac{\delta \bigl(\mathrm{TD}_2^u)^2 }{\delta V} \approx \frac{\delta \bigl(\mathrm{TD}_2^u)^2 }{\delta V}\bigl(X_t\bigr) = 2 \mathrm{TD}_2^u\bigl(X_t\bigr) \frac{\delta \mathrm{TD}_2^u}{\delta V} \bigl( X_t \bigr),
    \label{eq:TDgrad2}
\end{align}
where we evaluate the right-hand side term on a single realization of the trajectory to the ease the notation. In our numerical implementation, we use multiple trajectories to further improve the computation efficiency.
As we remark above, since \eqref{eqn:Ito2} holds true without taking the expectation, the right-hand side of \eqref{eq:TDgrad1} thus vanishes for the exact value function for any realization of $X_t$, in particular, the variance of the stochastic gradient is $0$. In comparison, while the stochastic estimate of \eqref{eq:TDgrad2} has the expectation $0$ for the exact value function, for each trajectory, the right-hand side is not $0$. This means that the stochastic gradient estimate \eqref{eq:TDgrad2} has a larger variance than the estimate \eqref{eq:TDgrad1}. Let us remark that the vanishing variance property of \eqref{eq:TDgrad1} is similar to quantum Monte Carlo~\cite{foulkes2001quantum}, for which the variance of the local energy estimate also vanishes at the ground state.  

In the following, to distinguish the two loss functions, we call the method based on $L_1$ \eqref{eqn:Critic1} the variance reduced LSTD (VR-LSTD), while that corresponding to $L_2$ \eqref{eqn:Critic2} is named the LSTD. We will demonstrate in our numerical experiments that VR-LSTD gives better results than LSTD.

\smallskip 

\subsubsection{Policy gradient for continuous optimal control problems}\label{sec:contPG}

For the actor part, we use policy gradient to improve the policy. According to the dynamical programming principle \cite{yong1999stochastic}, for the optimal value function $V$, we have
\begin{equation} \label{eqn:DP_cont}
V(X_0) = \inf_{u} \mathbb{E}_{u} \Bigl[ \int_0^{T} f(X_s,u(X_s)) e^{-\gamma s} \rd s + e^{-\gamma T} V(X_T) \mid X_0 \Bigr],
\end{equation}
where $u$ is minimized over the set of admissible controls. In other words, the control $u$ should minimize the functional on the right-hand side. Therefore, we can use the following loss function for the actor, for which we also incorporate the stopping time:
\begin{equation} \label{eqn:Actor}
J(u) = \mathbb{E}_{X_0 \sim \mu, u} \, \Bigl[\int_0^{T \land \tau} f(X_s,u(X_s)) e^{-\gamma s} \rd s +  \widehat{V}(X_{T \land \tau}) e^{-\gamma  (T \land \tau)}\Bigr],
\end{equation}
where $\widehat{V}$ is the current estimate of the value function (via TD learning in the critic part). Observe that this loss function is a continuous analog of \eqref{eq:MDP_actor}.

In the numerical algorithm, the control, as a high dimensional function, will be parametrized as a neural network $u(\cdot; \theta_u)$, where $\theta_u$ denotes collectively the parameters. The parameters are optimized using a stochastic approximation to gradients of $J(u)$.
Similarly to our discussion of policy gradient for the discrete case in Section~\ref{sec:MDP_PG}, when differentiating the loss function \eqref{eqn:Actor} w.r.t. the parameters of the control $\theta_u$, several terms would contribute to the derivative, including the control $u(\cdot)$ itself,  the SDE trajectory $X_{\cdot}$, the stopping time $\tau$, and also the estimated value function $\widehat{V}(\cdot)$. Similarly to the discrete case, we will drop the functional derivative of $\widehat{V}$ w.r.t.~$u$, i.e., the derivative $\frac{\delta \widehat{V}}{\delta u} \frac{\partial u}{\partial \theta_u}$, since the dependence of $\widehat{V}$ on $u$ is through the algorithm for the critic, e.g., the TD learning, which is impractical to track.
Furthermore, if $\hat{V}$ is the optimal value function, treating it as a fixed function and optimizing $u$ in \eqref{eqn:Actor} gives the optimal policy function.
Therefore, we approximate the functional derivative as 
\begin{multline}\label{eq:PG_explicit_gradient}
    \dfrac{\delta J}{\delta u}  \stackrel{\cdot}{=}   \mathbb{E}_{X_0 \sim \mu, u} \Biggl[ \int_0^{T \land \tau} \dfrac{\delta f(X_s, u(X_s)) }{\delta u} e^{-\gamma s} \rd s + \mathbbm{1}_{\{\tau < T\}} f(X_{\tau}, u(X_{\tau})) e^{-\gamma \tau} \dfrac{\delta \tau}{\delta u} \\
    + \nabla \widehat{V} (X_{T \land \tau}) e^{-\gamma  (T \land \tau)} \dfrac{\delta X_s}{\delta u}  \bigg\vert_{s = T \land \tau}     + \mathbbm{1}_{\{\tau < T\}} ( \mathcal{L}^u - \gamma) \widehat{V}(X_{\tau}) e^{-\gamma \tau} \dfrac{\delta \tau }{\delta u} \Biggr],
\end{multline}
where $\stackrel{\cdot}{=}$ indicates that we leave out the contribution from the functional derivative of $\widehat{V}$ w.r.t. $u$ and we have 
\begin{equation}
    \dfrac{\delta f(X_s, u(X_s)) }{\delta u} = \dfrac{\partial f}{\partial x} \dfrac{\delta X_s}{\delta u} + \dfrac{\partial f}{\partial u} \left( \mathrm{Id} + \nabla u(X_s)  \dfrac{\delta X_s}{\delta u} \right). 
\end{equation}
To obtain the formula, we have used It\^o's lemma to rewrite 
\begin{equation}\label{eq:Ito_exp}
    \mathbb{E}_{X_0 \sim \mu, u} [\widehat{V}(X_{T \land \tau}) e^{-\gamma(T \land\tau)}] = \mathbb{E}_{X_0 \sim \mu, u} \Bigl[\widehat{V}(X_0) + \int_0^{T\land\tau} ( \mathcal{L}^u - \gamma) \widehat{V}(X_s) \rd s\Bigr],
\end{equation}
and taken the derivative of the right-hand side w.r.t.{} $u$.

\section{Numerical algorithm}
In this section, we present our numerical algorithm for solving high dimensional HJB type elliptic PDEs based on the actor-critic framework discussed in the previous section. 

\subsection{Function approximation}
In order to numerically deal with the high dimensional functions $V$ and $u$, we use two neural networks to parametrize the value function $V(\cdot~; \theta_V)$ and the control $u(\cdot~; \theta_u)$, the parameters of which are denoted collectively by $\theta_V$ and $\theta_u$, respectively. We apply the structure of the residual neural network~\cite{he2016deep} in pursuit of better optimization performance. 
A neural network $\phi(x; \theta)$ with $l$ hidden layers is represented by
\begin{equation}
    \phi(x; \theta) = F_l \circ \sigma_l \circ F_{l-1} \circ \sigma_{l-1} \circ \cdots \circ F_1 \circ \sigma_1 \circ F_0 (x),
\end{equation}
where $F_i$ are linear transforms with dimensions depending on the width of hidden layers and the dimensions of inputs and outputs, and $\sigma_i$ are elementwise activate functions with skip connection: $\sigma_i(x) = x + \text{ReLU}(x)$. 

Moreover, note that the VR-LSTD loss function $L_1$ \eqref{eqn:Critic1} requires the gradient of the value function.
Since we are using a neural network parametrization  $V = V(\cdot~; \theta_V)$, a direct approach is to use autodifferentiation of $V(x;\theta_V)$ w.r.t. $x$ to calculate the gradient. 
We find that a better approach in practice is to use another neural network to represent $\nabla V$, which is consistent with the observations in \cite{han2018solving, han2020solving}. Thus, for VR-LSTD, the gradient of the value function is represented by a separate neural network $G(\cdot;\theta_{G})$ with collective parameters $\theta_{G}$. 

To summarize, w.r.t. the collective parameters, the loss functions for the critic corresponding to \eqref{eqn:Critic1} and \eqref{eqn:Critic2} are
\begin{align}
L_1(\theta_V, \theta_G)  & = \mathbb{E}_{X_0 \sim \mu, u} \,  \biggl[\Bigl( \int_0^{T \land \tau} e^{-\gamma s} f(X_s,u(X_s)) \rd s \label{eqn:Critic1_numerical}\\
& \hspace{10em} - \int_0^{T \land \tau} e^{-\gamma s} G(X_s; \theta_G)^{\top} \sigma(X_s,u(X_s)) \rd W_s \notag\\
& \hspace{10em} + e^{-\gamma (T \land \tau)}  V(X_{T \land \tau}; \theta_V) - V(X_0; \theta_V) \Bigr)^2\biggr],   \notag\\
L_2(\theta_V) & = \mathbb{E}_{X_0 \sim \mu, u} \,  \biggl[\Bigl(  \int_0^{T \land \tau} e^{-\gamma s} f(X_s,u(X_s)) \rd s + e^{-\gamma (T \land \tau)}  V(X_{T \land \tau}; \theta_V)\label{eqn:Critic2_numerical}\\
& \hspace{10em} - V(X_0; \theta_V) \Bigr)^2\biggr]. \notag
\end{align}
We remark that there is no need to add penalty terms in $L_1$ to ensure the consistency between $V(x,\theta_V)$ and $G(x;\theta_{G})$, because if we replace $V^u(\cdot)$ by $V(\cdot,\theta_V)$ in \eqref{eqn:Ito1} and plug it in \eqref{eqn:Critic1_numerical}, we have
\begin{align}
    L_1(\theta_V, \theta_G) & = \mathbb{E}_{X_0 \sim \mu, u} \,  \biggl[\Bigl( \int_0^{T \land \tau} e^{-\gamma s} [(\mathcal{L}^u V - \gamma V)(X_s; \theta_{V}) + f(X_s,u(X_s))] \rd s \label{eq:L1_consistency} \\
    & \hspace{2em} - \int_0^{T \land \tau} e^{-\gamma s} \left(\nabla_{x} V(X_s; \theta_{V}) - G(X_s; \theta_{G}) \right)^{\top} \sigma(X_s,u(X_s)) \rd W_s \Bigr)^2\biggr] \notag \\
    & = \mathbb{E}_{X_0 \sim \mu, u} \, \Bigl[\bigl( \int_0^{T \land \tau} e^{-\gamma s} [(\mathcal{L}^u V - \gamma V)(X_s; \theta_{V}) + f(X_s,u(X_s))] \rd s \bigr)^2\Bigr] \label{eq:term1}\\
    &\hspace{-0.5em} + \mathbb{E}_{X_0 \sim \mu, u} \, \Bigl[ \int_0^{T \land \tau} e^{-2\gamma s} \left|\sigma^{\top}(X_s,u(X_s)) \left(\nabla_{x} V(X_s; \theta_{V}) - G(X_s; \theta_{G}) \right)\right|^2 \rd s\Bigr],\label{eq:term2} 
\end{align}
where $\mathcal{L}^u$ is the generator of the SDE and we have used It\^o's isometry in the second step. Note that the first \eqref{eq:term1} and second \eqref{eq:term2} terms in \eqref{eq:L1_consistency} simultaneously enforce $V$ to be the value function and its gradient to be consistent with $G$.

For a neural network parametrization of $V$, it is not easy to directly impose the Dirichlet boundary condition $V = g$ on $\partial \Omega$ in the parametrization. Thus, instead, we add a penalty term to the loss functions \eqref{eqn:Critic1_numerical} or \eqref{eqn:Critic2_numerical} for the critic to help enforce the boundary condition 
\begin{equation}\label{eqn:boundary_loss}
\eta\, \mathbb{E}_{X \sim \text{Unif}(\partial \Omega)}\, \bigl[(V(X; \theta_V) - g(X))^2\bigr],
\end{equation}
where $\eta$ is a penalty hyperparameter and $\text{Unif}(\partial \Omega)$ denotes the uniform distribution on $\partial \Omega$.

\subsection{Discretization of SDEs and stochastic integrals}

In the implementation, we need to simulate numerically, based on a discretization of the diffusion process with approximating stopping time and exit point.
The solution to the PDE problem crucially depends on the boundary condition, and thus in control formulation, the exit time and position of the SDE at the boundary.
Several schemes have been developed in the literature to deal with the stopping time and exit point of the SDEs in related scenarios. Perhaps the most natural idea is to stop at the last step of the numerical SDE before exiting the domain, which has been tested in the context of using neural networks for solving PDEs in \cite{kremsner2020deep}. The error of such boundary treatment has been analyzed in \cite{gobet2000weak}. Moreover, several schemes have been proposed to improve the accuracy around the boundary.
In \cite{han2020derivative}, the authors approximate the exit position by the intersection of the domain boundary and the line segment between the consecutive two steps before and after exiting the domain. It has also been considered to reduce step size when the discretized trajectory approaches the boundary~\cite{buchmann2003solving}.
Some bias reduction schemes with the bubble-wrap or max-sampling exit condition are proposed in \cite{martin2021solving}.
After studying and testing several approaches for numerical discretization in our algorithms, we present two choices of discretization and give some remarks on the other schemes.

Let us start with a na\"ive approach. We can discretize the SDE\,\eqref{eqn:SDE} by the Euler--Maruyama scheme with a given partition of interval $[0,T]$: $0 = t_0 < t_1 < \cdots < t_N = T$, where a constant step size $\Delta t = \frac{T}{N}$ is used, so $t_n = n \Delta t$. The SDE is discretized as 
\begin{equation} \label{eqn:scheme_X}
\mathcal{X}_0 = X_0, ~~~~~ \mathcal{X}_{t_{n+1}} = \mathcal{X}_{t_n} + b(\mathcal{X}_{t_n},u_n) \Delta t + \sigma(\mathcal{X}_{t_n}, u_n) \xi_n \sqrt{\Delta t},
\end{equation}
where $u_n = u(\mathcal{X}_{t_n}; \theta_u)$ and $\xi_n \sim N(0, I_{d_w})$ follows the standard normal distribution.
Here, we use $\mathcal{X}_{t_n}$ to denote the discretized stochastic process, to distinguish from $X_t$, the continuous process. Given a numerical trajectory $\mathcal{X}_{t_n}, n = 0, \ldots, N$, we define 
\begin{equation}
    \bar{n} = \max \bigl\{ n \in \{0, \ldots, N\} \mid \mathcal{X}_{t_{i}} \in \Omega, ~ i = 0,1,\cdots,n
    \bigr\}. 
\end{equation}
Thus, if $\bar{n} < N$, $\mathcal{X}_{t_{\bar{n}+1}}$ exits the domain as $\mathcal{X}_{t_{\bar{n}+1}} \not \in \Omega$, while if $\bar{n} = N$ the trajectory $\mathcal{X}_{t_n}$ remains in the domain for $n = 0, 1, \ldots, N$. 

Perhaps the most direct and intuitive approach for the boundary treatment is to view $t = \bar{n}\Delta t$ as the stopping time, even though $\mathcal{X}_{t_{\bar{n}}}$ is still inside $\Omega$. 
This scheme will be referred to as the ``na\"ive scheme'' in the following. The stochastic integrations in \eqref{eqn:Critic1}, \eqref{eqn:Critic2}, and \eqref{eqn:Actor} are correspondingly approximated by
\begin{equation} \label{eqn:scheme_naive_integration}
\begin{aligned}
& \int_0^{T \land \tau} e^{-\gamma s} f(X_s,u_s) \rd s \approx \sum_{n=0}^{\bar{n}-1} e^{-\gamma n \Delta t} f(\mathcal{X}_{t_n}, u_n) \Delta t, \\
& \int_0^{T \land \tau} e^{-\gamma s} \nabla V(X_s)^{\top} \sigma(X_s, u_s) \rd W_s \approx \sum_{n=0}^{\bar{n}-1} e^{-\gamma n \Delta t} G(\mathcal{X}_{t_n}; \theta_{G})^{\top} \sigma(\mathcal{X}_{t_n}, u_n) \xi_n \sqrt{\Delta t}, 
\end{aligned}
\end{equation}
where $\xi_n \sqrt{\Delta t}$ is the same realization of Brownian increments as in \eqref{eqn:scheme_X}.
We remark that this discretization scheme is similar to the one used in \cite{kremsner2020deep} for solving degenerate semilinear elliptic equations, in particular, both algorithms approximate the stopping time by $\bar{n}\Delta t$. However, we aim to solve the value function in the whole domain,  while the method developed in \cite{kremsner2020deep} only aims at the value at a specific point; thus the overall framework of the algorithm is quite different. 

After discretization, the loss functions \eqref{eqn:Critic1}, \eqref{eqn:Critic2}, and \eqref{eqn:Actor} are further approximated by Monte Carlo samples: for each iteration, we draw $K$ independent sample trajectories ($K$ is known as the batch size) by drawing initial point $X_0$ from the distribution $\mu$ and independent increments of the Brownian motion. At each iteration, we also draw $K$ independent Monte Carlo samples uniformly from the boundary to approximate the expectation in \eqref{eqn:boundary_loss}.
To update the parameters of the neural networks, we employ the Adam optimizer \cite{Kingma2015adam}.

To apply the policy gradient method to the loss functional \eqref{eqn:Actor}, we need to differentiate the discretized functional w.r.t. the control, similar to the functional derivative setting considered above in  \eqref{eq:PG_explicit_gradient}. While the first and third terms in \eqref{eq:PG_explicit_gradient}, which involve derivatives of $J$ through its dependence on $u$ and the trajectory, can be easily dealt with on the discretized level using autodifferentiation, the second and fourth terms in \eqref{eq:PG_explicit_gradient} become tricky to deal with on the discrete level, since the stopping time is approximated by $\bar{n} \Delta t$, which is discrete so $\delta \bar{n} / \delta u$ is not really well defined. In our implementation, such terms are omitted in the policy gradient w.r.t. $u$; we leave a better numerical treatment of such terms to future works.

The pseudocode for our actor-critic method for solving high dimensional PDEs is summarized in Algorithm~\ref{alg1}.

\begin{algorithm}[ht]
\SetAlgoLined
\DontPrintSemicolon
\SetNoFillComment
 \SetKwInOut{Input}{input}\SetKwInOut{Output}{output}
\Input{A fully nonlinear PDE \eqref{eqn:NonlinearPDE}, terminal time $T$, number of time intervals $N$, loss weights $\eta$, neural network structures, number of iterations, learning rate, batch size $K$, the choice of TD}
\Output{Value function $V(\cdot~; \theta_V)$, its gradient $G(x;\theta_{G})$ if we choose \text{VR-LSTD}, and the control $u(\cdot~; \theta_u)$}
\medskip 
initialization: $\theta_{\text{value}}$ ($\theta_{\text{value}} = (\theta_V, \theta_G)$ for \text{VR-LSTD} and $\theta_{\text{value}} = \theta_V$ for \text{LSTD}) and $\theta_{u}$ \;
\medskip 
\For{\forcond}{
\tcc{critic steps}
Sample $K$ independent trajectories $\mathcal{X}^k_{t_n}, k = 1,2,\cdots,K$\;

Sample $K$ points on the boundary $\partial \Omega$ to enforce the boundary condition\;

Estimate the gradient of the chosen critic loss (\eqref{eqn:Critic1_numerical} $+$ \eqref{eqn:boundary_loss} or \eqref{eqn:Critic2_numerical} $+$ \eqref{eqn:boundary_loss}) w.r.t. $\theta_{\text{value}}$ using the $K$ trajectories and $K$ boundary points\;

Update parameters $\theta_{\text{value}}$ using the Adam optimizer\;

\tcc{actor steps}
Sample $K$ independent trajectories $\mathcal{X}^k_{t_n}, k = 1,2,\cdots,K$\;

Estimate the gradient of the actor loss \eqref{eqn:Actor} w.r.t. $\theta_{u}$ using the $K$ trajectories\;

Update parameters $\theta_{u}$ using the Adam optimizer\;
}
\caption{Neural network based actor-critic solver for fully nonlinear PDEs}
\label{alg1}
\end{algorithm}

\subsection{The adaptive step size scheme} 

It turns out in our numerical experiments that while the above na\"ive scheme is able to get reasonably accurate value functions, the approximation to control results in large errors, especially near the boundary (see Section~\ref{sec:numerical_results} for more details). To improve the accuracy near the boundary, we adaptively shrink the step size when the trajectory approaches the boundary $\partial \Omega$, instead of using the uniform time step size as in the na\"ive scheme. More specifically, we use the following scheme at the boundary, which is motivated by the integration scheme used in \cite{buchmann2003solving} for the Feynman--Kac representation of boundary value problems of the Poisson equation. 

The idea is to reduce the time step size adaptively when 
$\mathcal{X}_t$ is close to the boundary, and thus to improve the accuracy of the trajectory. We consider the  Euler--Maruyama scheme with varying step size given by 
\begin{equation}\label{eq:adaptive_scheme}
    \mathcal{X}_{t_{n+1}} = \mathcal{X}_{t_{n}} +   b(\mathcal{X}_{t_{n}},u_{n}) h(\mathcal{X}_{t_{n}}) + \sigma(\mathcal{X}_{t_{n}}, u_n) \sqrt{h(\mathcal{X}_{t_{n}})} \, \xi_{n},
\end{equation}
where the step size $h(\mathcal{X}_{t_n})$ depends on the current position of the trajectory. For the choice of step size, we define a subset near the boundary of $\Omega$ as
\begin{equation}\label{eq:gammaT}
    \Gamma = \bigl\{ x \in \overline\Omega ~|~ \text{dist}(x,\partial \Omega) \le \varsigma \sqrt{3d \Delta t} \bigr\},
\end{equation}
where $\varsigma = \sup_{x \in \Omega, u \in U} \lVert \sigma(x, u) \rVert$ is the supremum of the operator norm of $\sigma$. The adaptive choice of the step size is specified as follows:
\begin{enumerate}[wide,itemsep=1pt,topsep=1pt]
\item When $\mathcal{X}_{t_n} \in \Omega \backslash \Gamma$, it would be considered in the ``interior'' of $\Omega$, as it is very unlikely that after one time step with step size $\Delta t$ that the trajectory will exit the domain. Thus, we will use the basic constant step size $h(\mathcal{X}_{t_n}) = \Delta t$.
\item When $\mathcal{X}_{t_n} \in \Gamma$,  we decrease the step size according to the distance of the trajectory to the boundary, with the minimum step size set as $\frac{1}{10^4} \Delta t$:
\begin{equation*}
    h(\mathcal{X}_{t_n}) = \max \Bigl\{\frac{1}{3d \varsigma^2} \dist(\mathcal{X}_{t_n}, \partial \Omega)^2, \frac{1}{10^4} \Delta t \Bigr\}.
\end{equation*} 
This reduced step size, together with the width of $\Gamma$ defined in \eqref{eq:gammaT}, are decided such that
the probability that $\mathcal{X}_{t_n} \in \Omega \backslash \Gamma$ goes out of $\Omega$ in the next step is small. Note that when the step size is small, the diffusion dominates the drift term, and thus it suffices to incorporate the diffusion part in the choice. Note that we have used the supremum of $\lVert \sigma \rVert$ for simplicity, one could also choose the criteria more locally if $\sigma$ varies a lot across the domain. 
The minimum step size is set to balance the accuracy and computational cost as, otherwise, the scheme might spend an unnecessarily long time resolving the trajectory near the domain boundary.
\end{enumerate}

In summary, we choose the adaptive step size $h = h(\mathcal{X}_{t_n})$ as 
\begin{equation}\label{eqn:step_size}
h(\mathcal{X}_{t_{n}}) = 
\begin{cases} 
\Delta t, & \mathcal{X}_{t_n} \in \Omega \backslash \Gamma, \\
\max\{\frac{1}{3d \varsigma^2} \dist(\mathcal{X}_{t_n}, \partial \Omega)^2, \frac{1}{10^4} \Delta t\},  & \mathcal{X}_{t_n} \in \Gamma.
\end{cases}
\end{equation}
It should be noted that, as a result of the adaptive step size, different trajectories may have different discretized time steps.
The integrals are similarly discretized as in \eqref{eqn:scheme_naive_integration} with step size changed to $h(\mathcal{X}_{t_{n}})$: 
\begin{equation} \label{eqn:scheme_adaptive_integration}
\begin{aligned}
& \int_0^{T \land \tau} e^{-\gamma s} f(X_s,u) \rd s \approx \sum_{n=0}^{\bar{n}-1} e^{-\gamma \sum_{k=0}^{n-1} h(\mathcal{X}_{t_{k}})} f(\mathcal{X}_{t_n}, u_n) h(\mathcal{X}_{t_{n}}); \\
& \int_0^{T \land \tau} e^{-\gamma s} \nabla V(X_s)^{\top} \sigma(X_s) \rd W_s  \\
& \hspace{8em} \approx \sum_{n=0}^{\bar{n}-1} e^{-\gamma \sum_{k=0}^{n-1} h(\mathcal{X}_{t_{k}})} G(\mathcal{X}_{t_n}; \theta_{G})^{\top} \sigma(\mathcal{X}_{t_n},u_n) \xi_n \sqrt{h(\mathcal{X}_{t_{n}})}.
\end{aligned}
\end{equation}

For the policy gradient, similar to our numerical treatment in the case of na\"ive scheme,  we use autodifferentiation generated by the computational graph in practice, instead of directly numerically approximating the functional derivative defined in \eqref{eq:PG_explicit_gradient}. One reason is that the adaptive step size scheme further complicates the dependence of the trajectory and exit time on the control, compared with the na\"ive scheme and, hence, makes the direct numerical discretization of \eqref{eq:PG_explicit_gradient} even more difficult. In practice, the result from using autodifferentiation for the policy gradient seems to be quite accurate, as will be further discussed in the next section. 

\begin{remark}
In addition to the adaptive step size, in our numerical experiments, we have also tested the bounded sample of Brownian increments proposed in \cite{buchmann2003solving} to further avoid the potentially large error of the trajectory near the boundary due to tail events of the normal sample. We do not find, however, a significant difference in the result between using bounded samples versus the usual normal samples for Brownian increments. Therefore, we will stick to the normal samples for simplicity. 
\end{remark}

\begin{remark}
Moreover, besides adaptively shrinking the step size near the boundary,  we have tested two approaches using constant step size, but try to improve the estimate of the exit time and exit point of the na\"ive scheme instead. They do not yield satisfactory numerical results, so we will only briefly sketch the ideas without going into details or presenting numerical results. 

One scheme is adapted from \cite{han2020derivative}, which tries to determine the exit point on $\partial \Omega$ more accurately. 
In this scheme, the exit position $\mathcal{X}_{\tau}$ on $\partial \Omega$ is numerically approximated by the intersection of $\partial \Omega$ and the line segment between $\mathcal{X}_{t_{\bar{n}}}$ and $\mathcal{X}_{t_{\bar{n}+1}}$; the stopping time is correspondingly adjusted. 
The numerical result from the scheme is still not accurate enough for the control near the boundary. 

The linear interpolation above gives an error of order $\sqrt{\Delta t}$ due to the diffusion term; we can further improve the accuracy using a method proposed in  \cite{gobet2000weak}. Instead of linear interpolation, we seek for a coefficient $\rho \in (0,1]$ such that  $\mathcal{X}_{\tau}$, defined by
\begin{equation*}
\mathcal{X}_{\tau} = \mathcal{X}_{t_{\bar{n}}} +   b(\mathcal{X}_{t_{\bar{n}}},u_{\bar{n}}) \rho \Delta t + \sigma(\mathcal{X}_{t_{\bar{n}}},u_{\bar{n}}) \sqrt{\rho \Delta t} \, \xi_{\bar{n}},
\end{equation*}
is on $\partial \Omega$. 
In practice, we observe that numerically solving the coefficient $\rho$ makes the training unstable. 
\end{remark}

\section{Numerical examples} \label{sec:numerical_results}
In this section, we present the numerical results for the proposed method. We test on several examples: the linear quadratic regulator (LQR) problems, the stochastic Van der Pol oscillator problems, the diffusive Eikonal equations, and fully nonlinear elliptic PDEs derived from a regulator problem.
To test the performance of our algorithm, we do not assume knowledge of the true solution or the explicit formula for the control given the value function.
The considered dimensions in all four examples are as large as 20.
The algorithm is implemented in Python with the deep learning library TensorFlow 2.0~\cite{abadi2016tensorflow}.
In all the examples, the weight parameter $\eta$ associated with the boundary condition (cf.~\eqref{eqn:boundary_loss}) is set to 1 and the terminal time is $T=0.2$. The numbers of time intervals are $N=50$ for problems in $4$ dimensions ($4\mathrm{d}$) and $5\mathrm{d}$, and $N=100$ in $10\mathrm{d}$ and $20\mathrm{d}$. As for the architecture of the neural networks, the width of the hidden layers is set to 200 in all problems, while the numbers of hidden layers are 2 for problems in $4\mathrm{d}$ and $5\mathrm{d}$, and 3 in $10\mathrm{d}$ and $20\mathrm{d}$. 
During the training, we use piecewise constant learning rates of $\num{1e-3}$, $\num{1e-4}$, and $\num{1e-5}$ consecutively in order to achieve high accuracy. The numbers of steps with learning rate $\num{1e-3}$ are 20000 for problems in $4\mathrm{d}$, $5\mathrm{d}$, and $10\mathrm{d}$, and 30000 in $20\mathrm{d}$. The numbers of steps with learning rate $\num{1e-4}$ and $\num{1e-5}$ are both 10000 in the four examples. The batch sizes are $K=1024$ for problems in $4\mathrm{d}$ and $5\mathrm{d}$, and $K=2048$ in $10\mathrm{d}$ and $20\mathrm{d}$. The parameters in the numerical examples are determined empirically. In order to illustrate the effect of some parameters such as $T$ and the basic step size $\Delta t$, we also compare the results with different parameters in the first example.

During the training, we sample a validation set $\{X^k\}_{k=1}^K$ uniformly in $\Omega$, independent of the training, to evaluate the errors of the value function and the control.
Note that the validation size $K$ is the same as the batch size. We find that such sizes are enough to estimate the error accurately with a small variance.
The relative $L^2$ errors are computed by
\begin{equation} \label{eqn:err_value}
\text{err}_{V}^2 = \sum_{k=1}^K (V(X^k) - V(X^k; \theta_V))^2  / \sum_{k=1}^K V(X^k)^2
\end{equation}
and
\begin{equation} \label{eqn:err_control}
\text{err}_{u}^2 = \sum_{k=1}^K |u(X^k) - u(X^k; \theta_u)|^2  / \sum_{k=1}^K |u(X^k)|^2,
\end{equation}
where $V(\cdot)$ and $u(\cdot)$ are the true value and control functions, respectively (we will choose test examples such that these true solutions are known).
In addition to the errors above, we also visualize the density of the true value function and compare that with its neural network approximation, considering the difficulty of visualizing functions in high dimensions directly. Here, the density of a function $V$ is defined as the probability density function of $V(X)$, where $X$ is uniformly distributed in $\Omega$. In our numerical experiments, the density is estimated by Monte Carlo sampling. 

Our numerical results indicate that in all the examples, the value functions are approximated accurately, and the associated densities match well with that of the true solution. 
Furthermore, the numerical results show that, for the critic, the VR-LSTD performs better than LSTD, as expected. The adaptive step size scheme also significantly improves the accuracy, in particular, for the control.
The details can be found in the following subsections.
The code developed to solve these numerical examples is made publicly available on GitHub \cite{zhou2021actorcriticPDE-git}.

\subsection{LQR}
In this subsection we consider the PDE arising from the LQR problem, given by 
\begin{equation}
\Delta V(x) + \inf_{u \in \mathbb{R}^d} \bigl(\beta u^{\top} \nabla V(x) + p|x|^2 + q|u|^2 - 2kd \bigr) - \gamma V(x) =0 \quad \text{in } B_R \subset \mathbb{R}^d
\end{equation}
with boundary condition $V(x) = kR^2$ on $\partial B_R$, where $B_R = \{ x \in \mathbb{R}^d : |x| < R \}$.
Here  $p$, $q$, $\beta$, $k$ are positive constants such that
\begin{equation}
k = \frac{\sqrt{q^2 \gamma^2 + 4pq \beta^2} - \gamma q}{2\beta^2}.
\end{equation}
This is the HJB equation corresponding to the controlled stochastic process
\begin{equation}
\rd X_t = \beta u \rd t + \sqrt{2} \rd W_t
\end{equation}
with cost functional
\begin{equation}
J^u(x) = \mathbb{E} \Bigl[\int_0^{\tau} ( p|X_s|^2 + q|u(X_s)|^2 - 2kd ) e^{-\gamma s} \rd s + e^{-\gamma \tau} kR^2 \Bigr],
\end{equation}
where $\tau$ is the exit time of the domain $B_R$. The PDE has the exact solution as a quadratic function, $V(x) = k |x|^2$, and the optimal control is also explicitly given as $u^*(x) = \frac{-\beta}{2q} \nabla V(x) = \frac{- k \beta}{q} x$.

\begin{figure*}[t!]
    \centering
    \includegraphics[width=1.0\textwidth]{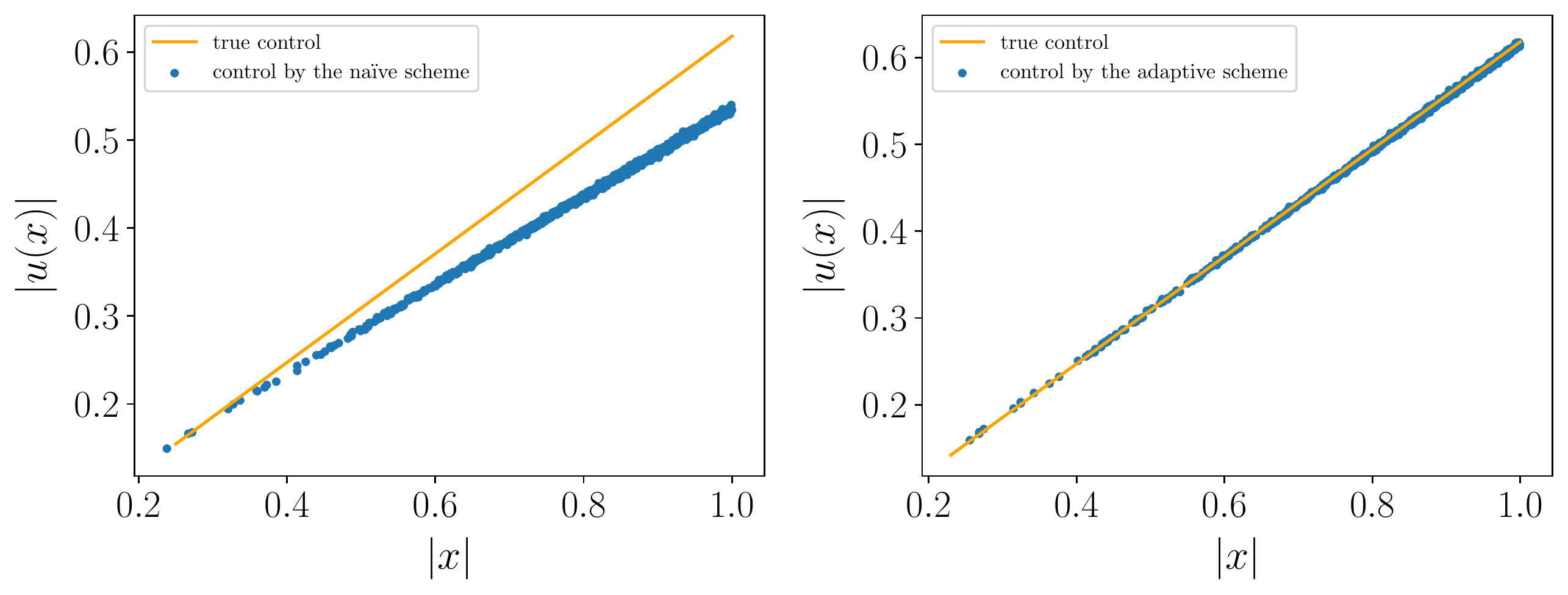}
    \caption{The comparison of two discretization schemes in the 5d LQR example. Left: the true optimal control and approximated control by the na\"ive scheme; right: the true optimal control and approximated control by the adaptive step size scheme. $x$-axis: the norm of $x$; $y$-axis: the norm of control $u$.}
    \label{fig:LQR_compare}
\end{figure*}

We choose the model parameters $p=q=R=\beta=\gamma=1$ and $k=(\sqrt{5}-1)/2$. The numerical results for our two versions of TDs with different discretization schemes in $5\mathrm{d}$ are shown in Table \ref{tab:error_lqr5d}.
The results of VR-LSTD have smaller errors due to the smaller asymptotic variance, as we discussed above. 
Moreover, the adaptive step size scheme is able to compute a more accurate control function, compared with the na\"ive scheme.
One possible reason is that the adaptive step size scheme samples more points $\mathcal{X}_{t_n}$ near the boundary, which helps to improve the accuracy of the control function near the boundary. To further illustrate the idea, let us compare the results for the na\"ive scheme and the adaptive step size scheme in $5\mathrm{d}$, both with critics optimized by VR-LSTD. 
The plot of the norm of the control $|u(x)|$ w.r.t.{} the norm of the variable $|x|$ is shown in Figure~\ref{fig:LQR_compare}. The error of the control for the na\"ive scheme is significantly larger near the boundary. The adaptive step size scheme achieves a uniform accuracy of the control in the whole domain.

\begin{figure*}[t!]
    \centering
    \includegraphics[width=0.98\textwidth]{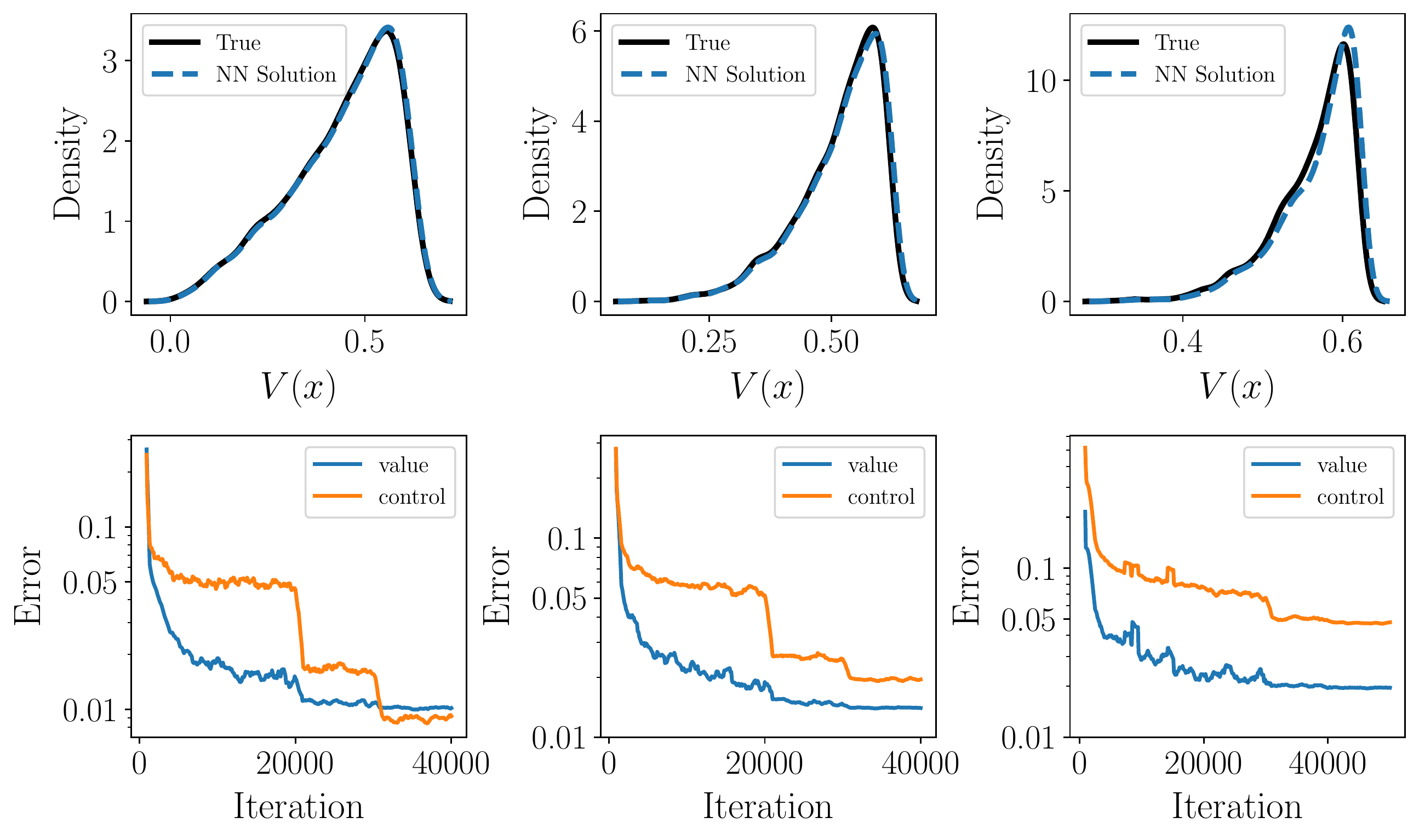}
    \caption{Top: density of $V$ for the LQR problem with $d=5$ (left), $d=10$ (middle), and $d=20$ (right).
    Bottom: associated error curves in the training process with $d=5$ (left), $d=10$ (middle), and $d=20$ (right).}
    \label{fig:lqr}
\end{figure*}

Therefore, for the rest of the numerical experiments, we will stick to the adaptive step size scheme and VR-LSTD loss function for the critic.
Figure~\ref{fig:lqr} shows the density and error curves for the LQR problem when $d=5, 10, 20$.
The sharp drop of the errors at steps 20000 and 30000 is due to the reduced learning rates at those steps.
The final errors of the value functions and controls are $\num{1.02e-2}$ and $\num{9.19e-3}$ in $5\mathrm{d}$; $\num{1.40e-2}$ and $\num{1.95e-2}$ in $10\mathrm{d}$; $\num{1.96e-2}$ and $\num{4.78e-2}$ in $20\mathrm{d}$. 

Considering that $T$ and the basic step size $\Delta t$ are two hyperparameters in our algorithm, we test different choices of their values in the $5\mathrm{d}$ LQR example and provide the errors in Table \ref{tab:compare_parameters}. The results show that our algorithm is not sensitive to the choice of $T$ and $\Delta t$.

\begin{table}[t!]
\centering
\begin{tabular}{c c | c c} 
\toprule
discretization  & TD variant & error of value function & error of control \\ 
\midrule
adaptive step size & VR-LSTD & {\textbf{\num{1.02e-2}}} & {\textbf{\num{9.19e-3}}} \\
adaptive step size & LSTD & \num{1.58e-1} & \num{1.17e-1}\\
na\"ive & VR-LSTD & \num{1.29e-2} & \num{1.24e-1} \\ 
na\"ive & LSTD & \num{1.41e-1} & \num{8.55e-2}\\
\bottomrule
\end{tabular}
\caption{Errors for different discretization schemes and TDs in $5\mathrm{d}$ LQR.}
\label{tab:error_lqr5d}
\end{table}

\begin{table}[!htbp]\centering 
\begin{tabular}{c|c|ccccc} 
\toprule
\multicolumn{2}{ c| }{$T$} & $T=0.04$ & $T=0.1$ & $T=0.2$ & $T=0.4$ & $T=0.8$ \\
\midrule
\multirow{2}{*}{50 time intervals}
&value & \num{8.40e-3} & \num{9.46e-3} & \num{1.04e-2} & \num{1.03e-2} & \num{9.97e-3}\\
&control & \num{1.15e-2} & \num{9.87e-3} & \num{1.01e-2} & \num{8.99e-3} & \num{8.33e-3}\\
\midrule 
\multirow{2}{*}{step size 0.004}
&value & \num{3.17e-2} & \num{1.40e-2} & \num{9.96e-3} & \num{9.16e-3} & \num{8.76e-3}\\
&control & \num{3.51e-2} & \num{1.19e-2} & \num{9.39e-3} & \num{1.03e-2} & \num{1.07e-2}\\
\bottomrule
\end{tabular}
\caption{
Errors of value and control functions with different parameters in the $5\mathrm{d}$ LQR example using adaptive step sizes and VR-LSTD. The first two rows denote different $T$, with the same number of time intervals $N=50$. The last two rows denote different $T$, with the same basic step size $\Delta t=0.004$. The relationship $T=N\Delta t$ always holds.
}
\label{tab:compare_parameters}
\end{table}

\subsection{Stochastic Van der Pol oscillator}
The Van der Pol oscillator is a popular example in the study of dynamical systems because of its chaotic behavior.
The stochastic Van der Pol oscillator has been studied in \cite{xu2011stochastic}, in which some internal or external noise is considered. In this subsection, we consider the generalized stochastic Van der Pol oscillator in high dimensional cases and solve the PDE
\begin{equation}
\Delta V(x) + \inf_{u \in \mathbb{R}^{d/2}} \bigl[ b(x,u)^{\top} \nabla V(x) + f(x,u)\bigr]  - \gamma V(x) = 0 \quad \text{in } B_R \subset \mathbb{R}^d,
\end{equation}
where $d=2n$ is even. The boundary condition is given by (with convention $x_0 = x_n$ and $x_{2n+1} = x_{n+1}$) \begin{equation}
    g(x) = a\sum_{i=1}^{2n} (x_{i})^2 - \epsilon \Bigl(\sum_{i=1}^n x_{i-1} x_i + \sum_{i=n+1}^{2n} x_i x_{i+1}\Bigr).
\end{equation}
Here $a$ and $\epsilon$ are positive constants. The drift field is given by
\begin{equation}
b_i(x,u) = \begin{cases}
x_{i+n} ~~ & (1 \le i \le n), \\
(1-x_{i-n}^2)x_i - x_{i-n} + u_{i-n} & (n+1 \le i \le 2n).
\end{cases}
\end{equation}
We choose the running cost as 
\begin{equation}\label{eq:runningcostvanderpol}
\begin{aligned}
f(x,u) =\;& q|u|^2 + \gamma [\sum_{i=1}^n (ax^2_i - \epsilon x_i x_{i-1}) + \sum_{i=n+1}^{2n} (ax^2_i - \epsilon x_i x_{i+1})] \\
&+ \frac{1}{4q} [(2ax_{n+1} - \epsilon x_{2n} - \epsilon x_{n+2})^2 + \sum_{i=n+2}^{2n} (2ax_i - \epsilon x_{i-1} - \epsilon x_{i+1})^2] - 4na\\
&- 2a \sum_{i=1}^n x_{n+i} x_i + \epsilon \sum_{i=1}^n x_{n+i} x_{i-1} + \epsilon \sum_{i=1}^{n-1} x_{n+i} x_{i+1} + \epsilon x_{2n} x_1\\
& - (x_{n+1} - x_1 - x_1^2 x_{n+1}) (2a x_{n+1} - \epsilon x_{2n} -\epsilon x_{n+2}) \\
&-\sum_{i=2}^n (x_{i+n} - x_i - x_i^2 x_{i+n}) (2a x_{i+n} - \epsilon x_{i+n-1} -\epsilon x_{i+n+1}),
\end{aligned}
\end{equation}
so that the true value function has an explicit formula:
\begin{equation}
V(x) = a\sum_{i=1}^{2n} (x_{i})^2 - \epsilon \Bigl(\sum_{i=1}^n x_{i-1} x_i + \sum_{i=n+1}^{2n} x_i x_{i+1}\Bigr).
\end{equation}
The corresponding optimal control is given by $u^*_1(x) =  - \frac{1}{2q} \partial_{n+1} V(x) = 2ax_{n+1} - \epsilon x_{2n} - \epsilon x_{n+2}$ and $u^*_i(x) = - \frac{1}{2q} \partial_{i+n} V(x) = 2ax_{i+n} - \epsilon x_{i+n-1} - \epsilon x_{i+n+1}$ for $i=2,3,\cdots,n$. 

The PDE can be reformulated as a stochastic control problem with the controlled SDE given by 
\begin{equation}
\rd X_t = b(X_t,u) \rd t + \sqrt{2} \rd W_t
\end{equation}
with objective function 
\begin{equation}
J^u(x) = \mathbb{E}\,\Bigl[\int_0^{\tau} f(X_s,u) e^{-\gamma s} \rd s + e^{-\gamma \tau} g(X_{\tau})\Bigr].
\end{equation}

\begin{figure*}[t!]
    \centering
    \includegraphics[width=0.98\textwidth]{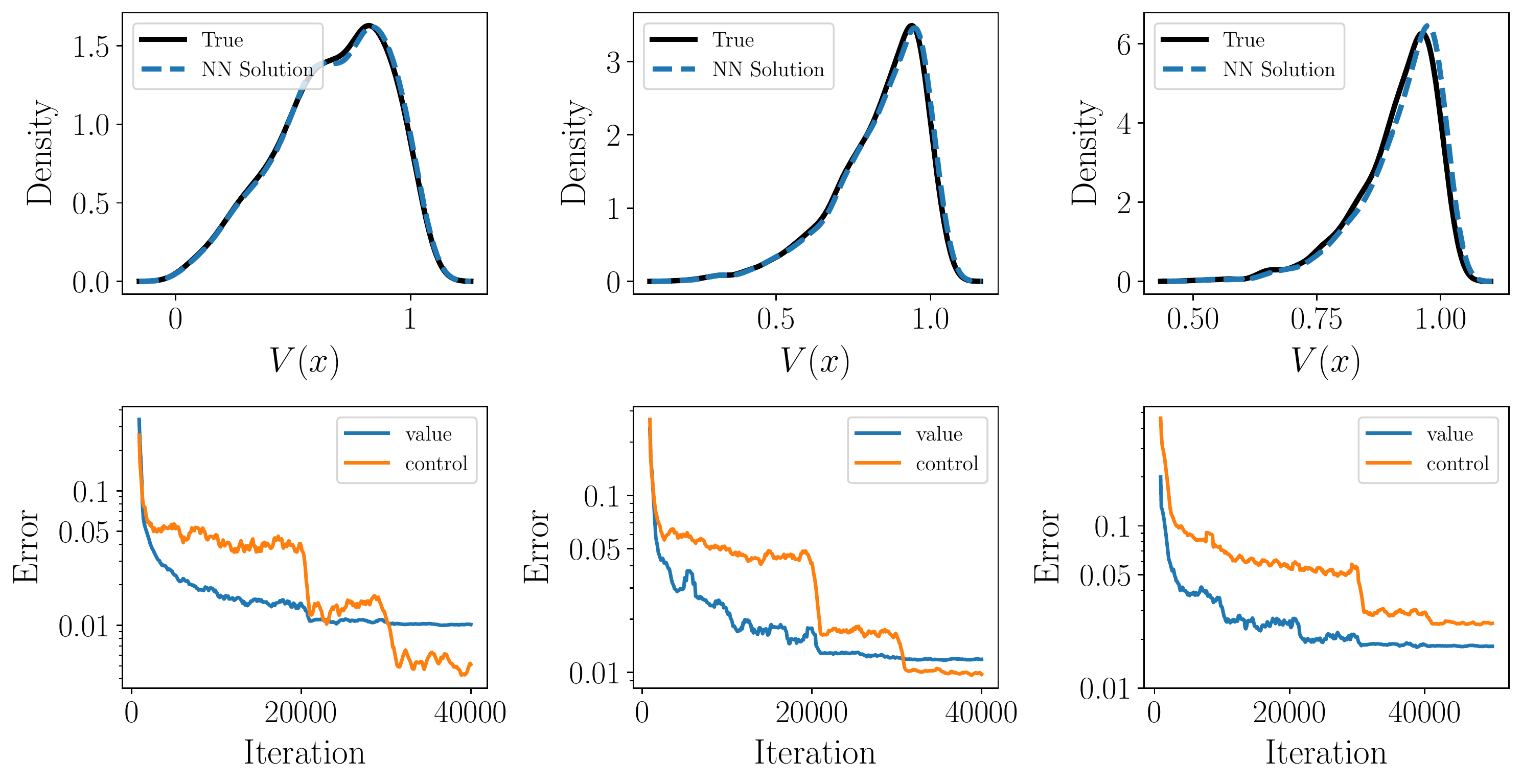}
    \caption{Top: density of $V$ for the Van der Pol problem with $d=4$ (left), $d=10$ (middle), and $d=20$ (right).
    Bottom: associated error curves in the training process with $d=4$ (left), $d=10$ (middle), and $d=20$ (right).}
    \label{fig:vdp}
\end{figure*}

In the numerical experiments, we take $a=q=R=\gamma=1$ and $\epsilon=0.1$.
Figure~\ref{fig:vdp} shows the density and error curves when $d=4, 10, 20$.
The algorithm learns reasonably nice shapes of the value functions.
The final errors of the value functions and controls are $\num{1.01e-2}$ and $\num{5.12e-3}$ in $4\mathrm{d}$; $\num{1.19e-2}$ and $\num{9.77e-3}$ in $10\mathrm{d}$; $\num{1.81e-2}$ and $\num{2.50e-2}$ in $20\mathrm{d}$.

\subsection{Diffusive Eikonal equation}
The Eikonal equation corresponds to the shortest-path problems with a given metric. In our experiments, we add a small diffusion term to regularize the equation (otherwise, the solution has kinks, which creates difficulty for the neural networks to approximate well in high dimensions). 
The diffusive Eikonal equation is given by 
\begin{equation}
\left\{ \begin{aligned}
& \epsilon \Delta V(x) + \inf_{u \in B_1} \bigl(c(x) u^{\top} \nabla V(x)\bigr) + 1  = 0 ~~~ \text{in} ~ B_R, \\
&V(x) = a_3-a_2 ~~~ \text{on} ~ \partial B_R,
\end{aligned} \right.
\end{equation}
where 
\begin{equation}
c(x) = \dfrac{3(d+1)a_3}{2d a_2(2a_2 - 3a_3|x|)} > 0
\end{equation}
is a real valued function. Here $a_2$ and $a_3$ are positive constants such that  $2a_2 - 3a_3R > 0$ and $\epsilon = 1/(2d a_2)$.
We choose the form of $c$ so that the true solution of the PDE is explicitly given by
\begin{equation}
V(x) = a_3|x|^3 - a_2|x|^2
\end{equation}
and the optimal control is $u^*(x) = x / |x|$.
In the numerical test, we take $a_2=1.2$, $a_3=0.2$, and $R=1$.

\begin{figure*}[t!]
    \centering
    \includegraphics[width=0.98\textwidth]{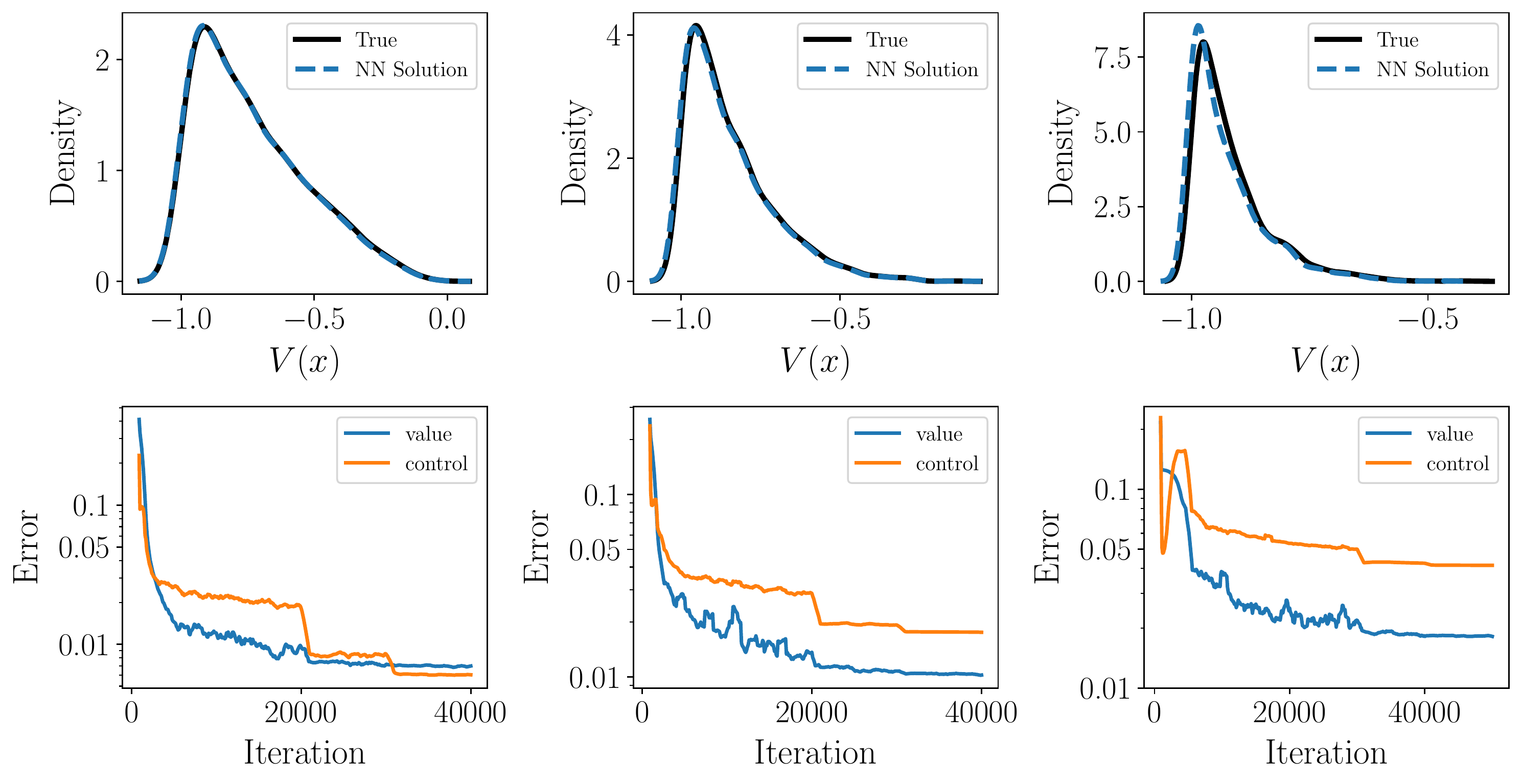}
    \caption{Top: density of $V$ for the Eikonal equation with $d=5$ (left), $d=10$ (middle), and $d=20$ (right).
    Bottom: associated error curves in the training process with $d=5$ (left), $d=10$ (middle), and $d=20$ (right).}
    \label{fig:ekn}
\end{figure*}

Unlike the previous two examples, the constraint on the control in this example poses a new challenge to the numerical algorithm. In order to ensure that the control $u$ is in the unit ball, we construct a specific structure of the neural network for the control. Instead of outputting the control directly, the neural network gives a $d + 1$ dimensional vector $(u_{\text{len}}, u_{\text{dir}}) \in \mathbb{R}^{d+1}$. The control is represented by
\begin{equation} \label{eqn:ekn_control}
u = \dfrac{u_{\text{dir}}}{\delta + \ReLU(u_{\text{len}}) + |u_{\text{dir}}|},
\end{equation}
where $\ReLU(x) = \max(0,x)$ and $\delta = 10^{-15}$.
This $\delta$ is to ensure that the denominator in \eqref{eqn:ekn_control} is not $0$ to prevent numerical singularity.
Figure~\ref{fig:ekn} shows the density and error curves for the Eikonal equation when $d=5, 10, 20$.
We also tried the straightforward parametrization of the control function as before, with an additional penalty term $\eta' \mathbb{E}_{X \sim \text{Unif}( \Omega)}[\text{ReLU}(|u(X)|-1)]$ in the loss for the actor. However, the numerical performances indicate that implementing the constraints of control directly like \eqref{eqn:ekn_control} is better than the penalty method.
The final errors of the value functions and controls are $\num{6.97e-3}$ and $\num{6.03e-3}$ in $5\mathrm{d}$; $\num{1.02e-2}$ and $\num{1.76e-2}$ in $10\mathrm{d}$; $\num{1.82e-2}$ and $\num{4.14e-2}$ in $20\mathrm{d}$.

\subsection{LQR with a nonconstant diffusion coefficient}
In this subsection, we consider a variant of the LQR in which the diffusion coefficient $\sigma$ is a function of both $x$ and $u$. Consider the HJB equation
\begin{equation}
\inf_{u \in \mathbb{R}^d} \bigl[\sum_{i=1}^d \left( \partial_i^2 V(x) (1+\epsilon x_i u_i)^2 + \beta \partial_i V(x) u_i \right) + q|u|^2 + \widetilde{f}(x) \bigr] - \gamma V(x) =0 \quad \text{in } B_R \subset \mathbb{R}^d,
\end{equation}
where
\begin{equation}
    \widetilde{f}(x) = \gamma k |x|^2 + \sum_{i=1}^d \dfrac{k^2 (\beta + 2 \epsilon)^2 x_i^2}{q + 2 k \epsilon^2 x_i^2} -2kd.
\end{equation}
In contrast to the previous three examples, this is a fully nonlinear PDE.
The corresponding SDE is 
\begin{equation}
    \rd X_t = \beta u_t \rd t + \sigma (X_t, u_t) \rd W_t,
\end{equation}
where $\sigma(x, u)$ is a diagonal matrix with $i$-th diagonal element  $\sqrt{2}(1 + \epsilon x_i u_i)$, $i=1,\cdots,d$. The running cost is $f(x, u) = q|u|^2 + \widetilde{f}(x)$. The true value function is $V(x) = k |x|^2$ and
the optimal control is 
\begin{equation}
    u_i^*(x) = -\dfrac{\beta \partial_i V(x) + 2 \epsilon x_i \partial_i^2 V(x)}{2q + 2 \epsilon^2 x_i^2 \partial_i^2 V(x)} = - \dfrac{(\beta + 2\epsilon) x_i}{q/k + 2\epsilon^2 x_i^2}.
\end{equation}
Note that this example coincides with the first example when $\epsilon =0$. In the numerical experiments, we set the parameters $q=R=\beta=\gamma=1, k=(\sqrt{5}-1)/2$ the same as the first example and $\epsilon=-1$. The final errors of the value functions and controls are $\num{9.98e-3}$ and $\num{1.63e-2}$ in $5\mathrm{d}$; $\num{1.50e-2}$ and $\num{4.95e-2}$ in $10\mathrm{d}$; $\num{1.96e-2}$ and $\num{5.25e-2}$ in $20\mathrm{d}$. This example showcases that our algorithm is able to solve fully nonlinear elliptic PDEs in high dimensions accurately.

\section{Conclusion and future directions}
\label{sec:conclusion}
In this paper, we propose and study numerical methods for high dimensional static HJB equations based on neural network parametrization and the actor-critic framework.
There are several promising directions for future research.
First, the scalability of the methods shall be further tested by problems of higher dimensions.
Second, it would be interesting to extend our methods to other types of boundary conditions like natural boundary conditions or broader types of equations, such as the porous medium equation. In both cases, the corresponding control formulation is not so clear. Third, one might explore the better numerical treatment of discretization of the functional derivative \eqref{eq:PG_explicit_gradient}, rather than relying on autodifferentiation. Finally, as an outstanding challenge in the field of deep learning, theoretical analysis for convergence and error analysis of the proposed numerical methods would be of great interest. 

\bibliographystyle{siamplain}
\bibliography{ref}

\begin{thebibliography}{10}

\bibitem{abadi2016tensorflow}
{\sc M.~Abadi, P.~Barham, J.~Chen, Z.~Chen, A.~Davis, J.~Dean, M.~Devin,
  S.~Ghemawat, G.~Irving, M.~Isard, et~al.}, {\em Tensorflow: A system for
  large-scale machine learning}, in 12th {USENIX} symposium on operating
  systems design and implementation ({OSDI} 16), 2016, pp.~265--283.

\bibitem{abdulla2007parametrized}
{\sc M.~S. Abdulla and S.~Bhatnagar}, {\em Parametrized actor-critic algorithms
  for finite-horizon {MDP}s}, in 2007 American Control Conference, IEEE, 2007,
  pp.~534--539.

\bibitem{barles2002convergence}
{\sc G.~Barles and E.~R. Jakobsen}, {\em On the convergence rate of
  approximation schemes for {H}amilton--{J}acobi--{B}ellman equations}, ESAIM:
  Mathematical Modelling and Numerical Analysis, 36 (2002), pp.~33--54.

\bibitem{beard1997galerkin}
{\sc R.~W. Beard, G.~N. Saridis, and J.~T. Wen}, {\em Galerkin approximations
  of the generalized {H}amilton--{J}acobi--{B}ellman equation}, Automatica, 33
  (1997), pp.~2159--2177.

\bibitem{beard1998approximate}
{\sc R.~W. Beard, G.~N. Saridis, and J.~T. Wen}, {\em Approximate solutions to
  the time-invariant {H}amilton--{J}acobi--{B}ellman equation}, Journal of
  Optimization theory and Applications, 96 (1998), pp.~589--626.

\bibitem{beck2020overview}
{\sc C.~Beck, M.~Hutzenthaler, A.~Jentzen, and B.~Kuckuck}, {\em An overview on
  deep learning-based approximation methods for partial differential
  equations}, arXiv preprint arXiv:2012.12348,  (2020).

\bibitem{becker2019deep}
{\sc S.~Becker, P.~Cheridito, and A.~Jentzen}, {\em Deep optimal stopping},
  Journal of Machine Learning Research, 20 (2019), p.~74.

\bibitem{bellman1966dynamic}
{\sc R.~Bellman}, {\em Dynamic programming}, Science, 153 (1966), pp.~34--37.

\bibitem{bhatnagar2006reinforcement}
{\sc S.~Bhatnagar and M.~S. Abdulla}, {\em A reinforcement learning based
  algorithm for finite horizon {M}arkov decision processes}, in Proceedings of
  the 45th IEEE Conference on Decision and Control, IEEE, 2006, pp.~5519--5524.

\bibitem{bhatnagar2009natural}
{\sc S.~Bhatnagar, R.~S. Sutton, M.~Ghavamzadeh, and M.~Lee}, {\em Natural
  actor-critic algorithms}, Automatica, 45 (2009), pp.~2471--2482.

\bibitem{boyan1999least}
{\sc J.~A. Boyan}, {\em Least-squares temporal difference learning}, in ICML,
  Citeseer, 1999, pp.~49--56.

\bibitem{buchmann2003solving}
{\sc F.~Buchmann and W.~Petersen}, {\em Solving {D}irichlet problems
  numerically using the {F}eynman--{K}ac representation}, BIT Numerical
  Mathematics, 43 (2003), pp.~519--540.

\bibitem{chan2019machine}
{\sc Q.~Chan-Wai-Nam, J.~Mikael, and X.~Warin}, {\em Machine learning for semi
  linear {PDE}s}, Journal of Scientific Computing, 79 (2019), pp.~1667--1712.

\bibitem{crandall1992user}
{\sc M.~G. Crandall, H.~Ishii, and P.-L. Lions}, {\em User’s guide to
  viscosity solutions of second order partial differential equations}, Bulletin
  of the American mathematical society, 27 (1992), pp.~1--67.

\bibitem{degris2012off}
{\sc T.~Degris, M.~White, and R.~Sutton}, {\em Off-policy actor-critic}, in
  International Conference on Machine Learning, 2012.

\bibitem{dolgov2019tensor}
{\sc S.~Dolgov, D.~Kalise, and K.~Kunisch}, {\em Tensor decompositions for
  high-dimensional {H}amilton--{J}acobi--{B}ellman equations}, arXiv preprint
  arXiv:1908.01533,  (2019).

\bibitem{duan2016benchmarking}
{\sc Y.~Duan, X.~Chen, R.~Houthooft, J.~Schulman, and P.~Abbeel}, {\em
  Benchmarking deep reinforcement learning for continuous control}, in
  International Conference on Machine Learning, PMLR, 2016, pp.~1329--1338.

\bibitem{han2016deep}
{\sc W.~E and J.~Han}, {\em Deep learning approximation for stochastic control
  problems}, arXiv preprint arXiv:1611.07422,  (2016).

\bibitem{weinan2017deep}
{\sc W.~E, J.~Han, and A.~Jentzen}, {\em Deep learning-based numerical methods
  for high-dimensional parabolic partial differential equations and backward
  stochastic differential equations}, Communications in Mathematics and
  Statistics, 5 (2017), pp.~349--380.

\bibitem{han2020algorithms}
{\sc W.~E, J.~Han, and A.~Jentzen}, {\em Algorithms for solving high
  dimensional {PDE}s: From nonlinear monte carlo to machine learning}, arXiv
  preprint arXiv:2008.13333,  (2020).

\bibitem{forsyth2007numerical}
{\sc P.~A. Forsyth and G.~Labahn}, {\em Numerical methods for controlled
  {H}amilton--{J}acobi--{B}ellman {PDE}s in finance}, Journal of Computational
  Finance, 11 (2007), p.~1.

\bibitem{foulkes2001quantum}
{\sc W.~Foulkes, L.~Mitas, R.~Needs, and G.~Rajagopal}, {\em Quantum {M}onte
  {C}arlo simulations of solids}, Reviews of Modern Physics, 73 (2001), p.~33.

\bibitem{gobet2000weak}
{\sc E.~Gobet}, {\em Weak approximation of killed diffusion using euler
  schemes}, Stochastic processes and their applications, 87 (2000),
  pp.~167--197.

\bibitem{han2020deep}
{\sc J.~Han and R.~Hu}, {\em Deep fictitious play for finding markovian {N}ash
  equilibrium in multi-agent games}, in Mathematical and Scientific Machine
  Learning, PMLR, 2020, pp.~221--245.

\bibitem{han2021recurrent}
{\sc J.~Han and R.~Hu}, {\em Recurrent neural networks for stochastic control
  problems with delay}, arXiv preprint arXiv:2101.01385,  (2021).

\bibitem{han2020convergence}
{\sc J.~Han, R.~Hu, and J.~Long}, {\em Convergence of deep fictitious play for
  stochastic differential games}, arXiv preprint arXiv:2008.05519,  (2020).

\bibitem{han2018solving}
{\sc J.~Han, A.~Jentzen, and W.~E}, {\em Solving high-dimensional partial
  differential equations using deep learning}, Proceedings of the National
  Academy of Sciences, 115 (2018), pp.~8505--8510.

\bibitem{han2020convergenceBSDE}
{\sc J.~Han and J.~Long}, {\em Convergence of the deep bsde method for coupled
  {FBSDE}s}, Probability, Uncertainty and Quantitative Risk, 5 (2020),
  pp.~1--33.

\bibitem{han2020solving}
{\sc J.~Han, J.~Lu, and M.~Zhou}, {\em Solving high-dimensional eigenvalue
  problems using deep neural networks: A diffusion monte carlo like approach},
  Journal of Computational Physics, 423 (2020), p.~109792.

\bibitem{han2020derivative}
{\sc J.~Han, M.~Nica, and A.~R. Stinchcombe}, {\em A derivative-free method for
  solving elliptic partial differential equations with deep neural networks},
  Journal of Computational Physics, 419 (2020), p.~109672.

\bibitem{he2016deep}
{\sc K.~He, X.~Zhang, S.~Ren, and J.~Sun}, {\em Deep residual learning for
  image recognition}, in Proceedings of the IEEE conference on computer vision
  and pattern recognition, 2016, pp.~770--778.

\bibitem{henry2017deep}
{\sc P.~Henry-Labordere}, {\em Deep primal-dual algorithm for {BSDE}s:
  Applications of machine learning to {CVA} and {IM}}, Available at SSRN
  3071506,  (2017).

\bibitem{ji2020deep}
{\sc S.~Ji, S.~Peng, Y.~Peng, and X.~Zhang}, {\em Deep learning method for
  solving stochastic optimal control problem via stochastic maximum principle},
  arXiv preprint arXiv:2007.02227,  (2020).

\bibitem{ji2020three}
{\sc S.~Ji, S.~Peng, Y.~Peng, and X.~Zhang}, {\em Three algorithms for solving
  high-dimensional fully coupled {FBSDE}s through deep learning}, IEEE
  Intelligent Systems, 35 (2020), pp.~71--84.

\bibitem{kalise2018polynomial}
{\sc D.~Kalise and K.~Kunisch}, {\em Polynomial approximation of
  high-dimensional {H}amilton--{J}acobi--{B}ellman equations and applications
  to feedback control of semilinear parabolic {PDE}s}, SIAM Journal on
  Scientific Computing, 40 (2018), pp.~A629--A652.

\bibitem{kang2017mitigating}
{\sc W.~Kang and L.~C. Wilcox}, {\em Mitigating the curse of dimensionality:
  sparse grid characteristics method for optimal feedback control and {HJB}
  equations}, Computational Optimization and Applications, 68 (2017),
  pp.~289--315.

\bibitem{Kingma2015adam}
{\sc D.~P. Kingma and J.~Ba}, {\em Adam: {A} method for stochastic
  optimization}, in 3rd International Conference on Learning Representations,
  {ICLR} 2015, San Diego, CA, USA, May 7-9, 2015, Conference Track Proceedings,
  Y.~Bengio and Y.~LeCun, eds., 2015, \url{http://arxiv.org/abs/1412.6980}.

\bibitem{klebaner2005introduction}
{\sc F.~C. Klebaner}, {\em Introduction to stochastic calculus with
  applications}, World Scientific Publishing Company, 2005.

\bibitem{konda2000actor}
{\sc V.~R. Konda and J.~N. Tsitsiklis}, {\em Actor-critic algorithms}, in
  Advances in Neural Information Processing Systems, Citeseer, 2000,
  pp.~1008--1014.

\bibitem{kremsner2020deep}
{\sc S.~Kremsner, A.~Steinicke, and M.~Sz{\"o}lgyenyi}, {\em A deep neural
  network algorithm for semilinear elliptic {PDE}s with applications in
  insurance mathematics}, Risks, 8 (2020), p.~136.

\bibitem{kunisch2004hjb}
{\sc K.~Kunisch, S.~Volkwein, and L.~Xie}, {\em {HJB}-{POD}-based feedback
  design for the optimal control of evolution problems}, SIAM Journal on
  Applied Dynamical Systems, 3 (2004), pp.~701--722.

\bibitem{lygeros2004reachability}
{\sc J.~Lygeros}, {\em On reachability and minimum cost optimal control},
  Automatica, 40 (2004), pp.~917--927.

\bibitem{maei2010toward}
{\sc H.~R. Maei, C.~Szepesv{\'a}ri, S.~Bhatnagar, and R.~S. Sutton}, {\em
  Toward off-policy learning control with function approximation}, in
  Proceedings of the 27th International Conference on Machine Learning
  (ICML-10), 2010, pp.~719--726.

\bibitem{martin2021solving}
{\sc C.~Martin, H.~Zhang, J.~Costacurta, M.~Nica, and A.~R. Stinchcombe}, {\em
  Solving elliptic equations with brownian motion: Bias reduction and temporal
  difference learning}, Methodology and Computing in Applied Probability,
  (2021), pp.~1--24.

\bibitem{mitchell2005time}
{\sc I.~M. Mitchell, A.~M. Bayen, and C.~J. Tomlin}, {\em A time-dependent
  {H}amilton--{J}acobi formulation of reachable sets for continuous dynamic
  games}, IEEE Transactions on automatic control, 50 (2005), pp.~947--957.

\bibitem{mitchell2003overapproximating}
{\sc I.~M. Mitchell and C.~J. Tomlin}, {\em Overapproximating reachable sets by
  {H}amilton--{J}acobi projections}, Journal of Scientific Computing, 19
  (2003), pp.~323--346.

\bibitem{mnih2013playing}
{\sc V.~Mnih, K.~Kavukcuoglu, D.~Silver, A.~Graves, I.~Antonoglou, D.~Wierstra,
  and M.~Riedmiller}, {\em Playing atari with deep reinforcement learning},
  arXiv preprint arXiv:1312.5602,  (2013).

\bibitem{nakamura2019adaptive}
{\sc T.~Nakamura-Zimmerer, Q.~Gong, and W.~Kang}, {\em Adaptive deep learning
  for high-dimensional {H}amilton--{J}acobi-{B}ellman equations}, arXiv
  preprint arXiv:1907.05317,  (2019).

\bibitem{nusken2020solving}
{\sc N.~N{\"u}sken and L.~Richter}, {\em Solving high-dimensional
  {H}amilton--{J}acobi--{B}ellman {PDE}s using neural networks: perspectives
  from the theory of controlled diffusions and measures on path space}, arXiv
  preprint arXiv:2005.05409,  (2020).

\bibitem{osher1988fronts}
{\sc S.~Osher and J.~A. Sethian}, {\em Fronts propagating with
  curvature-dependent speed: Algorithms based on {H}amilton--{J}acobi
  formulations}, Journal of computational physics, 79 (1988), pp.~12--49.

\bibitem{oster2019approximating}
{\sc M.~Oster, L.~Sallandt, and R.~Schneider}, {\em Approximating the
  stationary {H}amilton--{J}acobi--{B}ellman equation by hierarchical tensor
  products}, arXiv preprint arXiv:1911.00279,  (2019).

\bibitem{pardoux1998backward}
{\sc {\'E}.~Pardoux}, {\em Backward stochastic differential equations and
  viscosity solutions of systems of semilinear parabolic and elliptic {PDE}s of
  second order}, in Stochastic Analysis and Related Topics VI, Springer, 1998,
  pp.~79--127.

\bibitem{pardoux1990adapted}
{\sc E.~Pardoux and S.~Peng}, {\em Adapted solution of a backward stochastic
  differential equation}, Systems \& Control Letters, 14 (1990), pp.~55--61.

\bibitem{pardoux1992backward}
{\sc E.~Pardoux and S.~Peng}, {\em Backward stochastic differential equations
  and quasilinear parabolic partial differential equations}, in Stochastic
  partial differential equations and their applications, Springer, 1992,
  pp.~200--217.

\bibitem{peng1991probabilistic}
{\sc S.~Peng}, {\em Probabilistic interpretation for systems of quasilinear
  parabolic partial differential equations}, Stochastics and Stochastics
  Reports, 37 (1991), pp.~61--74.

\bibitem{pereira2019learning}
{\sc M.~A. Pereira, Z.~Wang, I.~Exarchos, and E.~A. Theodorou}, {\em Learning
  deep stochastic optimal control policies using forward-backward {SDE}s}, in
  Robotics: science and systems, 2019.

\bibitem{peters2008natural}
{\sc J.~Peters and S.~Schaal}, {\em Natural actor-critic}, Neurocomputing, 71
  (2008), pp.~1180--1190.

\bibitem{pham2021neural}
{\sc H.~Pham, X.~Warin, and M.~Germain}, {\em Neural networks-based backward
  scheme for fully nonlinear {PDE}s}, SN Partial Differential Equations and
  Applications, 2 (2021), pp.~1--24.

\bibitem{rao2009survey}
{\sc A.~V. Rao}, {\em A survey of numerical methods for optimal control},
  Advances in the Astronautical Sciences, 135 (2009), pp.~497--528.

\bibitem{richardson2006numerical}
{\sc S.~Richardson and S.~Wang}, {\em Numerical solution of
  {H}amilton--{J}acobi--{B}ellman equations by an exponentially fitted finite
  volume method}, Optimization, 55 (2006), pp.~121--140.

\bibitem{silver2016mastering}
{\sc D.~Silver, A.~Huang, C.~J. Maddison, A.~Guez, L.~Sifre, G.~Van
  Den~Driessche, J.~Schrittwieser, I.~Antonoglou, V.~Panneershelvam, and
  M.~Lanctot}, {\em Mastering the game of {Go} with deep neural networks and
  tree search}, Nature, 529 (2016), pp.~484--489.

\bibitem{silver2014deterministic}
{\sc D.~Silver, G.~Lever, N.~Heess, T.~Degris, D.~Wierstra, and M.~Riedmiller},
  {\em Deterministic policy gradient algorithms}, in International Conference
  on Machine Learning, PMLR, 2014, pp.~387--395.

\bibitem{sutton2018reinforcement}
{\sc R.~S. Sutton and A.~G. Barto}, {\em Reinforcement learning: An
  introduction}, MIT press, 2018.

\bibitem{vamvoudakis2010online}
{\sc K.~G. Vamvoudakis and F.~L. Lewis}, {\em Online actor-critic algorithm to
  solve the continuous-time infinite horizon optimal control problem},
  Automatica, 46 (2010), pp.~878--888.

\bibitem{wang2003numerical}
{\sc S.~Wang, L.~S. Jennings, and K.~L. Teo}, {\em Numerical solution of
  {H}amilton--{J}acobi--{B}ellman equations by an upwind finite volume method},
  Journal of Global Optimization, 27 (2003), pp.~177--192.

\bibitem{wang2016sample}
{\sc Z.~Wang, V.~Bapst, N.~Heess, V.~Mnih, R.~Munos, K.~Kavukcuoglu, and
  N.~de~Freitas}, {\em Sample efficient actor-critic with experience replay},
  arXiv preprint arXiv:1611.01224,  (2016).

\bibitem{xu2011stochastic}
{\sc Y.~Xu, R.~Gu, H.~Zhang, W.~Xu, and J.~Duan}, {\em Stochastic bifurcations
  in a bistable {D}uffing--{V}an der {P}ol oscillator with colored noise},
  Physical Review E, 83 (2011), p.~056215.

\bibitem{yong1999stochastic}
{\sc J.~Yong and X.~Zhou}, {\em Stochastic controls: Hamiltonian systems and
  {HJB} equations}, vol.~43, Springer, 1999.

\bibitem{zhou2021actorcriticPDE-git}
{\sc M.~Zhou, J.~Han, and J.~Lu}, {\em Actor-critic method for high dimensional
  static {H}amilton--{J}acobi--{B}ellman partial differential equations based
  on neural networks}.
\newblock \url{https://github.com/MoZhou1995/DeepPDE\_ActorCritic}.

\end{thebibliography}
\end{document}